\documentclass[final,12pt]{article}
\usepackage{verbatim}
\usepackage{latexsym}
\usepackage{amsmath,amsthm,amsfonts}
\usepackage{eucal}
\usepackage{cases}

\usepackage{mathtext}
\usepackage[cp1251]{inputenc}
\usepackage[T2A]{fontenc}
\usepackage[dvips]{graphicx}
\usepackage{amsmath}
\usepackage{amssymb}
\usepackage{amsxtra}
\usepackage{latexsym}
\usepackage{ifthen}
\usepackage{color}

\def\beq{\begin{equation}}
\def\eeq{\end{equation}}
\def\bea{\begin{eqnarray}}
\def\ena{\end{eqnarray}}

\long\def\drop#1{}

\newtheorem{theorem}{Theorem}
\newtheorem{definition}{Definition}[section]

\numberwithin{equation}{section}

\begin{document}

\title{\textbf{Use of mathematical modeling to study pressure regimes in normal and Fontan blood flow circulations}}

\author{J.P. Keener$^1$, M. Chugunova$^2$, R.M. Taranets$^3$, M.G. Doyle$^4$\\
 {\small $^1$Department of Mathematics, University of Utah,} \\ {\small  $^2$Institute of Mathematical Sciences, Claremont Graduate University,} \\
 {\small $^3$Institute of Applied Mathematics and Mechanics} \\
 {\small of National Academy of Sciences of Ukraine,} \\
 {\small $^4$ Department of Mechanical and Industrial Engineering,} \\
 {\small and Division of Vascular Surgery, University of Toronto.}}
\maketitle

\begin{abstract}
We develop two mathematical lumped parameter models for blood pressure distribution in the Fontan blood flow circulation:
an ODE based spatially homogeneous model and a PDE based spatially inhomogeneous model. We present numerical simulations
of the cardiac pressure-volume cycle and study the effect of pulmonary resistance on cardiac output.
We analyze solutions of two initial-boundary value problems for a non-linear parabolic partial differential equation (PDE models)
with switching in the time dynamic boundary conditions which model blood pressure distribution in the cardiovascular system
with and without Fontan surgery. We also obtain necessary conditions for parameter values of the PDE models for existence and uniqueness of non-negative bounded periodic solutions.
\end{abstract}

\section{Introduction to Fontan surgical procedure}
In a normal biventricular heart, the systemic and pulmonary blood circulations are in series and each one is supported
by a ventricle. The Fontan surgical procedure is applied to a malformed heart that is characterized by the presence of only one (left or right) functional ventricular chamber, and for which a biventricular repair is not possible \cite{mondesert2013fontan}. The Fontan procedure was first introduced in 1968 and involved
routing systemic venous blood flow to the pulmonary arteries, without use of an additional ventricular pumping
chamber. In its modern form, the Fontan procedure involves creating an anastomosis (surgical connection) between the superior vena cava (the vein responsible for carrying blood from the upper body to the heart) and the left and right pulmonary arteries, and between the inferior vena cava (the vein responsible for carrying blood from the lower body to the heart) and the left and right pulmonary arteries, with the connection between the inferior vena cava and the pulmonary arteries being made using a synthetic graft called an extracardiac conduit. The resulting connection between these four vessels is referred to as the total cavopulmonary connection. The Fontan circulation creates an unusual state in which the force driving the pulmonary blood flow is solely a residue (in the
systemic venous pressure) of the only functional ventricular chamber's contractile force. This single ventricle propels blood flow
through the systemic arteries and capillaries, with the systemic venous return passively entering the pulmonary circulation.
The ventricle is doing nearly twice the expected amount of work because it has to pump blood to both the body and the lungs.

Fontan surgery is an extraordinary story of success in that it has allowed a generation of newborn
babies with the most severe forms of congenital heart disease to survive into adulthood (estimated
prevalence of approximately $1$ per $3000$ births) \cite{khairy2007univentricular}. Though life-saving,
a univentricular Fontan circulation does not, however, reproduce biventricular physiology and has
been considered abnormal in the sense that systemic venous hypertension (mean pressure > 10 mm Hg) occurs
simultaneously with pulmonary arterial hypotension (mean pressure < 15 mm Hg) \cite{de2005fontan}.
In patients with Fontan physiology, life expectancy remains far below projected age- and sex-matched
normative values. Patients with Fontan procedures most frequently die from heart failure or from thromboemboli \cite{khairy2008long}.
The incidence of thromboembolic deaths rises sharply $15$ years after Fontan surgery \cite{valente2009congenital, d2007fontan}.
At the same time heart failure deaths are very uncommon during the first $10$ years after Fontan surgery, with a steady decline
in survival thereafter. Associated factors include protein-losing enteropathy, single right
ventricle morphology, and increased Fontan pressures \cite{mondesert2013fontan}.

Late Fontan failure might progress gradually over years with an absence of overt symptoms. Fontan patients have lived with less than ideal cardiac
output their entire lives and might not recognize decline in functional status until deterioration is significantly advanced. In the medical literature,
failure of the Fontan circulation is divided into $3$ main categories: ventricular dysfunction, systemic complications of Fontan physiology,
and chronic Fontan failure \cite{goldberg2011failing}. In a cross-sectional analysis of 546 children with Fontan
procedures, $27\%$ had abnormal ventricular ejection fractions and $72\%$ had diastolic dysfunction. The prevalence of
systolic and diastolic ventricular dysfunction continues to increase in adulthood \cite{piran2002heart, eicken2003hearts}.
Over the past $5-10$ years, a number of studies have described the effect of Fontan physiology on the liver. Hepatic
venous pressure after Fontan surgery might be $3-4$ times higher than normal, with levels commensurate with congestive
heart failure in adults. Implications of this elevated hepatic venous pressure over the long term remain to be fully
understood \cite{kendall2008hepatic, camposilvan2008liver}. Fontan physiology is characterized by progressively
decreasing cardiac output and increasing central venous pressure over time. The average peak oxygen consumption ranges
from $19$ to $28$ mL/kg per minute, or $50\%-60\%$ of predicted values \cite{fernandes2010serial, paridon2008cross}.
In the third decade of life, hospitalization rates and symptoms increase significantly \cite{diller2010predictors}.
In the future, a "Fontan pump" or a "Fontan assist device" might usher in a new era of ventricular replacement.
Although investigations have begun, such devices still appear to remain many years away \cite{rodefeld2010cavopulmonary, bhavsar2009intravascular}.
Recently, the use of the right ventricular Impella has been described to direct flow into the
pulmonary artery, resulting in a modest reduction in central venous pressure \cite{haggerty2012experimental}.

\section{Previous mathematical modeling results for Fontan circulation}
\subsection{Computational fluid dynamics}
Computational fluid dynamics is a powerful tool that can be used to gain insight into the local blood flow dynamics in the Fontan circulation. These simulations are used to model the detailed 3-D hemodynamics of a particular region in the body, such as the total cavopulmonary connection, rather than the complete circulation.  A simplified three-dimensional model was used \cite{taylor1996computational} to simulate the local fluid dynamics for different designs of the total cavopulmonary connection, allowing a quantitative evaluation  of the dissipated energy in each of the examined configurations. The authors show that, from a comparative point of view,
the energetic losses can be greatly reduced if a proper hydraulic design of the connection is adopted, which also allows control of the blood flow distribution into the lungs.

Under the assumptions that vessel walls are completely rigid (according to surgical reports, the vessel diameter
change per cardiac cycle is around $5-10 \%$ in most of the major arteries) and all vessels are symmetric, numerical simulations of blood flow to the lung after a surgical Fontan procedure are described in \cite{dubini1996numerical}.
The authors simulated the full nonlinear Navier-Stokes equations using a streamline finite element method. They analyzed the blood
flow dynamics for different values for the offset between the superior vena cava and the inferior vena cava anastomoses and concluded that the optimal distance for the offset is about $7$ mm.
It was shown in \cite{bazilevs2009computational} that wall flexibility can play an important role in determining quantities of hemodynamic interest in the Fontan connection.  However, \cite{tang2017} recently showed that while fluid-structure interaction effects are important for instantaneous quantities of interest in the Fontan circulation, they have a negligible impact on time-averaged values.

According to the article \cite{degroff2008modeling}, the main quantities of importance in modeling the Fontan procedure
are:\begin{itemize}
      \item Vessel diameters and flow rates representative of the
range seen in the patient group under study including
resting and exercise states
      \item Vessel sizes and flow rates matched appropriately
      \item Compliant vessels, accurate modeling of surgical
anastomosis sites, and surgical material used (unless
proven unnecessary)
      \item Unsteady flow
      \item Effects of respiration
      \item Correctly shaped vessel anatomy
    \end{itemize}

Two different types of boundary conditions, time-averaged and pulsatile, were analyzed in \cite{wei2016can}.
The authors derive a patient-specific sensitivity criterion which provides a guideline for determining when
time-averaged boundary conditions can be used to save computational time.

Recent advances in imaging methods and patient-specific modeling now reveal increasingly detailed information about
blood flow patterns in healthy and diseased Fontan states. Building on these tools, there is now an opportunity to couple blood flow
simulation with optimization algorithms to improve the design of surgeries and devices, incorporating more information
about the flow physics in the design process to augment current medical knowledge. To do so, there is a need for efficient optimization tools that are appropriate for unsteady fluid mechanics problems, particularly for
the optimization of complex patient-specific models in the presence of uncertainty. The state of the art in optimization
tools for virtual surgery, device design, and model parameter identification in cardiovascular flow and mechanobiology applications
are reviewed in \cite{yang2010constrained}. In this work, the authors perform optimization on a model Y-graft design problem.
This work represents the first use of formal design optimization methods for the Fontan surgery, and also demonstrates the
applicability of the optimization framework on a pulsatile flow problem with multiple design parameters and constraints.

\subsection{Lumped parameter models}
While computational fluid dynamics models can be used to calculate detailed three-dimensional blood flow in the total cavopulmonary connection, the computational costs of this approach prevent it from being used to simulate the entire circulatory system.  Because Fontan failure is a systemic problem, reduced order methods, such as lumped parameter models, can be used to study this phenomenon.  Lumped parameter models are based on the analogy between fluid flow and electric circuits.  In addition to modelling the entire circulatory system, these models can also be used to generate appropriate upstream and downstream boundary conditions for computational fluid dynamics simulations.

A lumped parameter model of the Fontan circulation was used by \cite{throckmorton2011} to generate boundary conditions for a computational fluid dynamics model used to design a Fontan assist device. In a study by \cite{liang2014}, lumped parameter models of the Fontan circulation and the normal circulation were compared to determine differences between the two circulations in the regulations of cardiac output and central venous pressures. In a study by \cite{kung2014}, a lumped parameter model was used to study the Fontan circulation under exercise conditions.

The objective of the present study is to develop lumped parameter models of the Fontan circulation with the goal of understanding the systematic changes that occur during Fontan failure. The outline of this paper is as follows.  Section 3 describes an ODE model of the Fontan circulation and presents some basic results from this model, showing its consistency with physiological behaviour.  Section 4 describes PDE modeling of the normal and  Fontan circulations, which is an extension of the ODE approach that allows for spatially inhomogeneous variation of the model parameters.  Section 5 determines range values of parameters for PDE models from the Section 3 such that a unique non-nenegative solution exists. We also show in the appendix how to construct super- and sub-solutions for PDE models.

\section{Spatially homogeneous ODE model of blood pressure distribution for the Fontan circulation}
\begin{figure}
\centerline{\includegraphics[height=6.0cm]{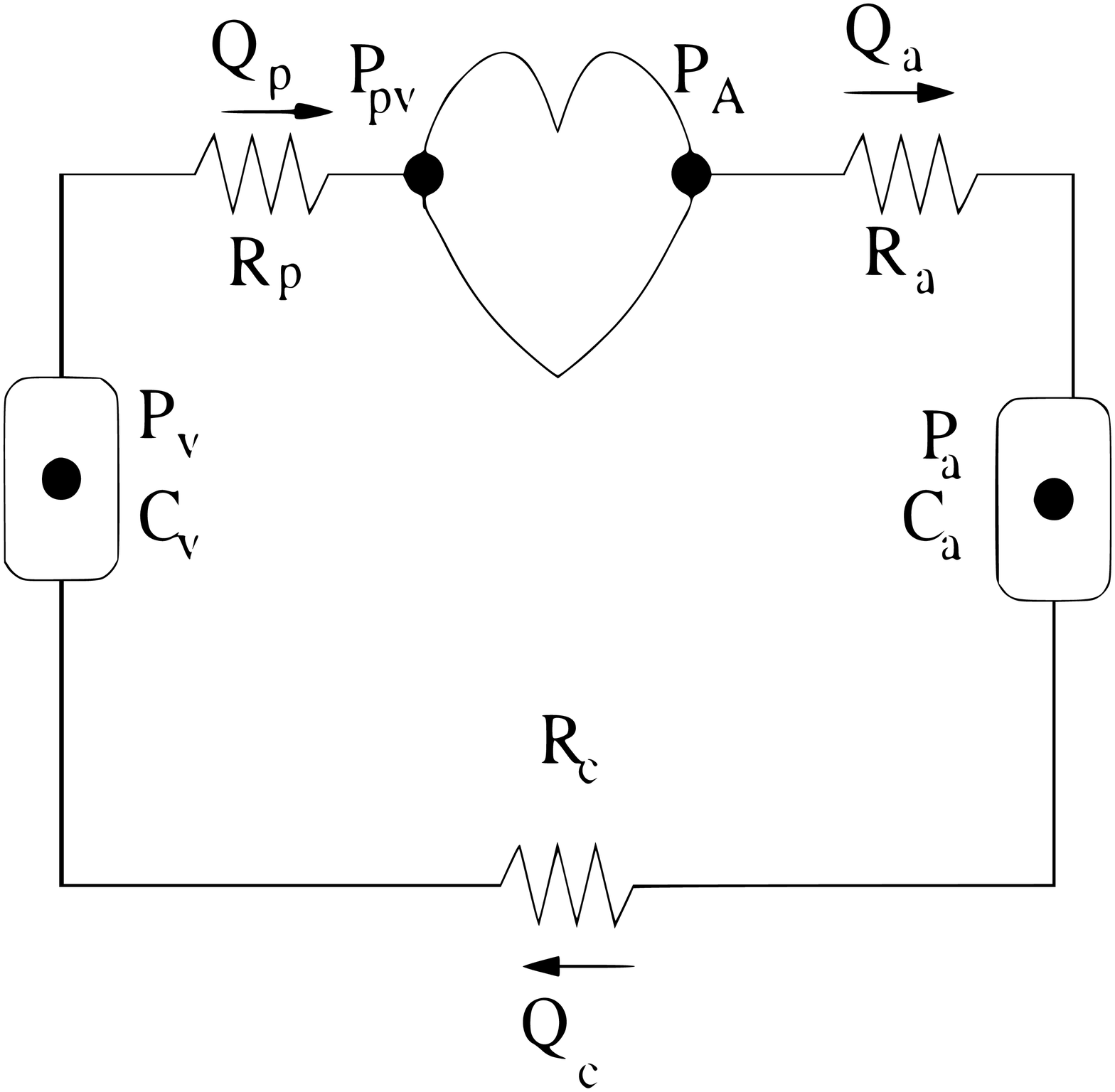}\hspace{1cm}\includegraphics[height=5.5cm]{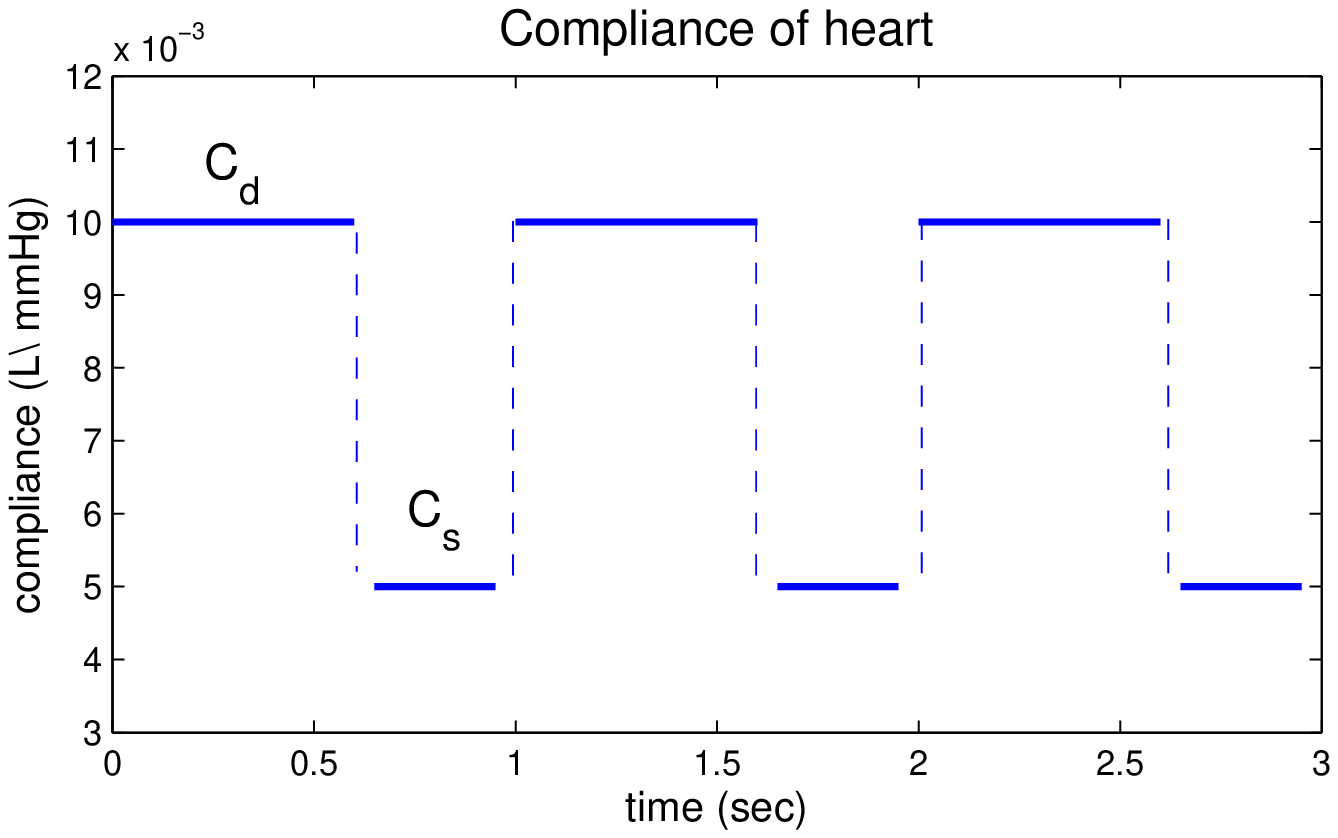}}
\caption{Circuit diagram for the Fontan heart circulatory system (left) and piecewise heart compliance as a function of time (right).}\label{fig:fig1}
\end{figure}
A simple model of the Fontan circulation can be based on an electric circuit approach.
This model consists of five compartments: the heart, the arterial system, the capillary system, the venous
system, and the pulmonary system (lungs). For the Fontan circulation, all compartments are connected in series around a single loop. In particular the flow must pass sequentially through the arteries, capillaries, veins and lungs, before returning to the heart, as shown in Fig.~\ref{fig:fig1}.

We model the capillary and pulmonary systems as linear resistance vessels. That is,
we assume that the pressure drop across the vessel is proportional to the flow through the
vessel, with a constant of proportionality called the resistance, labeled  $R_c$  and $R_p$, for the
capillary and pulmonary systems, respectively. We assume these vessels have no compliance.
We model the arterial and venous systems as linear compliance vessels, in which the volume
of the vessel is proportional to the pressure in the vessel, with constant of proportionality
called the compliance. We allow for compliance vessels to have resistance, which is modeled
in the same way as for the resistance vessels. The compliances of the arterial and venous
systems are labeled $C_a$ and $C_v$, respectively, while their resistances are $R_a$ and $R_v$. The
heart is considered to be a linear compliance vessel with different compliances depending on
whether it is relaxed (in diastole) or contracted (in systole). As shown in Fig.~\ref{fig:fig1}, the heart has compliance $C_d$ in diastole ($0\leq t < 0.7$) and compliance $C_s$ in systole ($0.7\leq t <1$). We assume that all vessels contain some basal volume at zero pressure, which are denoted as $V_i^0$.

\begin{figure}
\centerline{ \includegraphics[height=6cm]{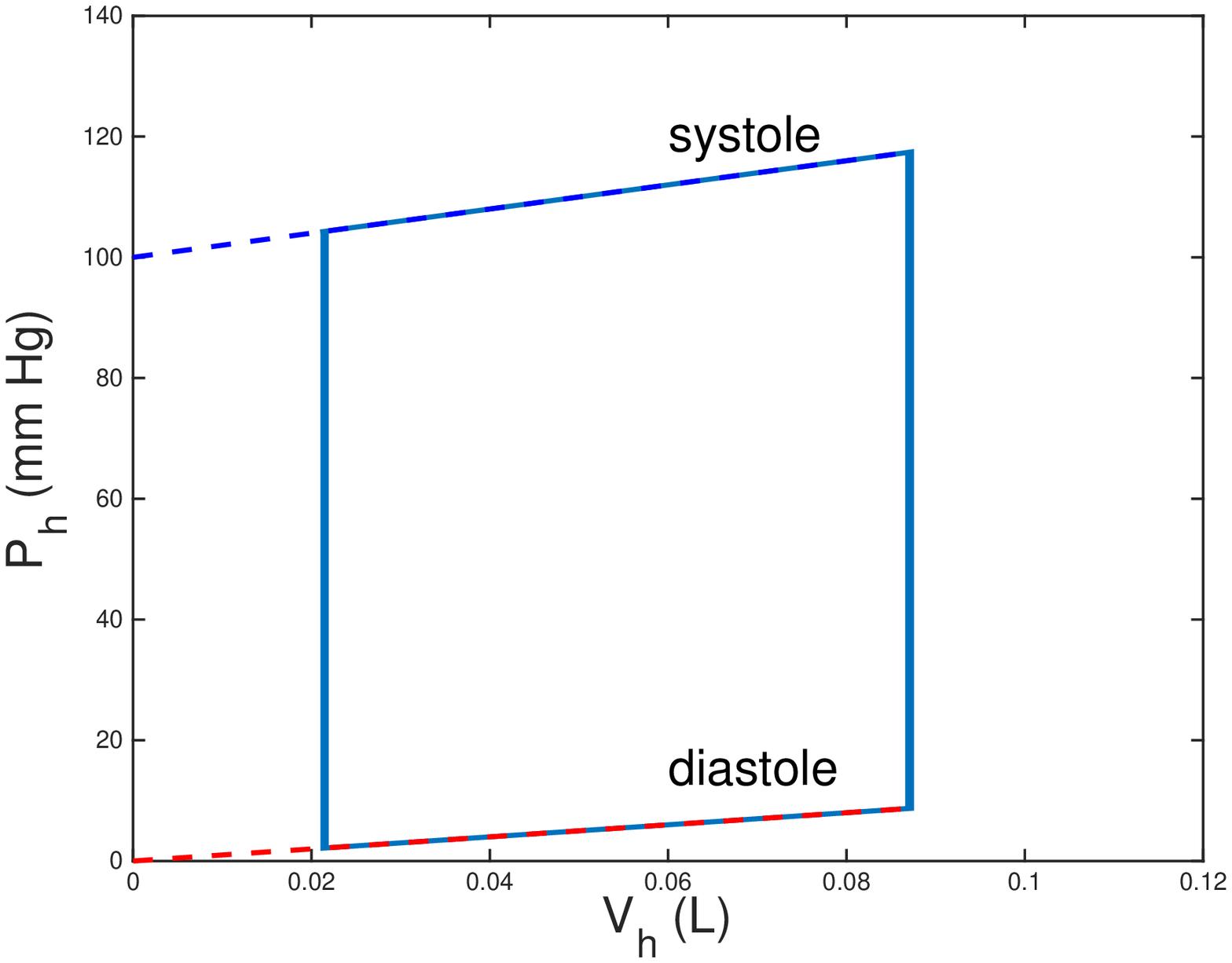}\hspace{1cm}\includegraphics[height=6cm]{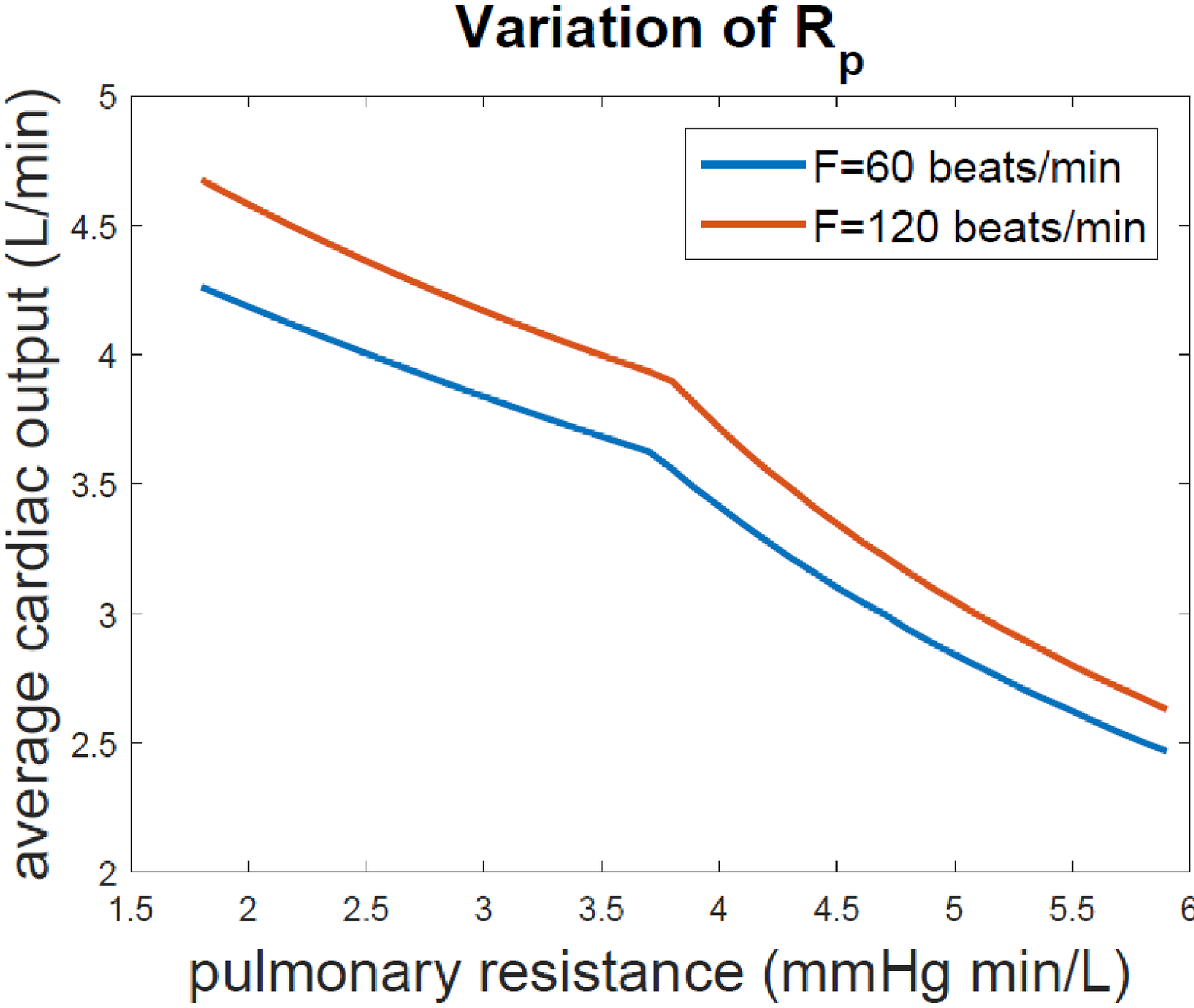}}
\caption{The cardiac pressure-volume cycle from our model (left) and effect of pulmonary resistance on cardiac output(right).}\label{fig:fig2}
\end{figure}

The variables in the system are the volumes $V_a$,
$V_v$ and $V_h$ of the compliance vessels (the arterial system, the venous system and the heart), and the
pressures $P_A$, $P_a$, $P_v$ and $P_{pv}$ at different points along the loop; see Fig.~\ref{fig:fig1}. The parameters of the system are the resistances, compliances and basal volumes of all the vessels as well as the total blood volume $V_T$. Estimates
for all the parameters can be found by measurement on individuals.
For the results that follow, the values of the parameters that we use have been taken from the literature and are given in Table~\ref{tab:tab1}.
Parameter values:  $V$ in units of L, $C$ in units of L/mm Hg and $R$ in units of mm Hg$\cdot$ min/L. It should be noted that $V_s$ was chosen to be negative to achieve a reasonable value of systolic heart compliance $C_s$.  In practice, the systolic pressure is never a small value or $0$, so the volume in the heart never reaches the negative value $V_s = -0.5$ L.
\begin{center}
\begin{table}[h]
\caption {Parameter values} \label{tab:tab1}
\hspace{0.5cm}
\begin{center}
\begin{tabular}{|c|c|c|c|}
\hline
$V_a^0$&1.0&$C_a$&0.00125\\
\hline
$V_v^0$&2.5& $C_v$&0.0625\\
\hline
$V_d^0$&0&$C_d$&0.01\\
\hline
$V_s^0$&-0.5&$C_s$&0.005\\
\hline
$R_a$&0.5&$R_c$&17\\
\hline
$R_p$&2.9&$V_T$& 5.0\\
\hline
\end{tabular}
\end{center}

\end{table}
\end{center}
We consider the compliance of the heart to be a piecewise
constant function, with value $C_s$ in systole and value $C_d$ in diastole. To ensure appropriate directionality of the forcing, we assume that there are "perfect" valves where the pulmonary vein enters the heart and where the aorta leaves the heart. Anatomically, in Fontan circulation, the pulmonary vein enters a common atrium which is separated from the single ventricle by an atrioventricular valve.  Depending on whether the patient has a functioning right or left ventricle, this valve is either the tricuspid valve or the mitral valve.  Similarly, the single ventricle is separated from the aorta by a semilunar valve, which is either the pulmonary valve or the aortic valve. For the purposes of these models, we are neglecting the common atrium by considering it an extension of the pulmonary vein, and we refer to the valve separating the pulmonary vein from the heart as the "pulmonary vein valve" and the valve separating the heart from the aorta as the "aortic valve". In our model, during diastole, the pulmonary vein valve is open and the aortic valve is closed, allowing flow into the heart, and during systole, the pulmonary vein valve is closed while the aortic valve is open allowing flow out of the heart and into the arteries. We assume that the valves are perfect in the sense that they open and close instantaneously and in synchrony, and that they do not allow back flow, regardless of pressure differences or flow characteristics. These assumptions lead to a cardiac pressure-volume cycle as presented in  Fig.~\ref{fig:fig2}.

We use conservation of volume in each compartment to set up the dynamic equations. In particular, we must have the rate of change in volume of a compartment equal to the difference of the flow in and flow out. Due to the assumptions on the valves, there are discontinuities in the variables and their derivatives as the valves
switch. Thus, it is convenient to distinguish the time intervals in which the heart is in systole
and in which it is in diastole. As such, there is one set of model equations that is valid
in diastole, while another set of equations is valid in systole.

The three flow rates are the arterial, capillary and pulmonary flow rates, defined as
\begin{equation}
Q_a=\tfrac{P_A-P_a}{R_a},\quad
Q_c=\tfrac{P_a-P_v}{R_c},\quad
Q_p=\tfrac{P_v-P_{pv}}{R_p}.
\end{equation}
The three compliance volumes are the arterial, venous and heart volumes, defined as
\begin{equation}
V_a = V_a^0+C_aP_a,\quad
V_v = V_v^0+C_vP_v,\quad
V_h = V_h^0+C_hP_h,
\end{equation}
with
\beq
V_h = \left\{\begin{array}{cc}V_d^0+ C_dP_{pv}&{\rm during ~ diastole}\\V_s^0+C_sP_A&{\rm during~ systole}\end{array}\right..
\eeq

Consequently, the dynamic equations are, during inflow (diastole),

\begin{equation}
\left\{
  \begin{array}{ll}\vspace{0.2cm}
   C_a\tfrac{dP_a}{dt} &=  Q_a-Q_c, \\ \vspace{0.2cm}
   C_v\tfrac{dP_v}{dt}&= Q_c-Q_p, \\ \vspace{0.2cm}
   C_d\tfrac{dP_{pv}}{dt} &=Q_p,
  \end{array}
\right.
\end{equation}
so that $Q_a=0$ and $P_A = P_a$,
and during outflow (systole)
\begin{equation}
\left\{
  \begin{array}{ll}\vspace{0.2cm}
   C_a\tfrac{dP_a}{dt} &=  Q_a-Q_c, \\ \vspace{0.2cm}
   C_v\tfrac{dP_v}{dt}&= Q_c-Q_p,\\ \vspace{0.2cm}
   C_s\tfrac{P_A}{dt} &= -Q_a,
  \end{array}
\right.
\end{equation}
so that $Q_p = 0$ and $P_v = P_{pv}$.  In addition, total volume must always be conserved, so that
\beq \label{eq:conv}
V_T = V_a+V_v+V_h.
\eeq

It is possible to simplify these equations by using the conservation law (\ref{eq:conv}).  In particular,
during inflow (diastole),
\begin{equation}
\left\{
  \begin{array}{ll}\vspace{0.2cm}
   C_a\tfrac{dP_a}{dt} &=  -\tfrac{P_a-P_v}{R_c}, \\\vspace{0.2cm}
   C_v\tfrac{dP_v}{dt}&= \tfrac{P_a-P_v}{R_c}-\tfrac{P_v-P_{pv}}{R_p},\\\vspace{0.2cm}
   P_{pv} &= \tfrac{1}{C_d}(V_T-V_a^0-C_aP_a-V_v^0-C_vP_v-V_d^0),\\\vspace{0.2cm}
   P_A &= P_a,
  \end{array}
\right.
\end{equation}
and during outflow (systole)
\begin{equation}
\left\{
  \begin{array}{ll}\vspace{0.2cm}
   C_a\tfrac{dP_a}{dt} &= \tfrac{P_A-P_a}{R_a}-\tfrac{P_a-P_v}{R_c},\\\vspace{0.2cm}
   C_v\tfrac{dP_v}{dt}&=\tfrac{P_a-P_v}{R_c},\\\vspace{0.2cm}
   P_A&=\tfrac{1}{C_s}(V_T-V_a^0-C_aP_a-V_v^0-C_vP_v-V_s^0),\\\vspace{0.2cm}
   P_{pv} &= P_v.
  \end{array}
\right.
\end{equation}

In addition to diastole and systole, the cardiac cycle consists of two isovolumetric phases, during which both heart valves are closed, and the heart undergoes a change in pressure in response to a change in its shape, while maintaining a constant blood volume.  Isovolumetric contraction occurs following diastole, during which the heart muscle contracts, increasing the pressure until the aortic valve opens to start systole.  Isovolumetric relaxation occurs following systole, during which the heart muscle relaxes, decreasing the pressure until the pulmonary vein valve opens.  In our model, which instantaneously switches from diastole to systole, the heart pressure $P_h$ is determined in a way that is consistent with the isovolumetric constraint.  In particular, during systole,
\begin{equation}
\left\{
  \begin{array}{ll}\vspace{0.2cm}
   V_h = V_T -(V_a^0+V_v^0+C_aP_a+C_vP_v),\\\vspace{0.2cm}
   P_h = (V_h-V_s^0)/C_s,
  \end{array}
\right.
\end{equation}
and during diastole,
\begin{equation}
\left\{
  \begin{array}{ll}\vspace{0.2cm}
   V_h = V_T - (V_a^0+V_v^0+C_aP_a+C_vP_v),\\\vspace{0.2cm}
   P_h = (V_h-V_d^0)/C_d.
  \end{array}
\right.
\end{equation}
Consequently, while $P_a$ and $P_v$ are continuous functions of time, $P_h$ experiences jump discontinuities at the transitions between diastole and systole.

\begin{figure}
\centerline{ \includegraphics[height=5.4cm]{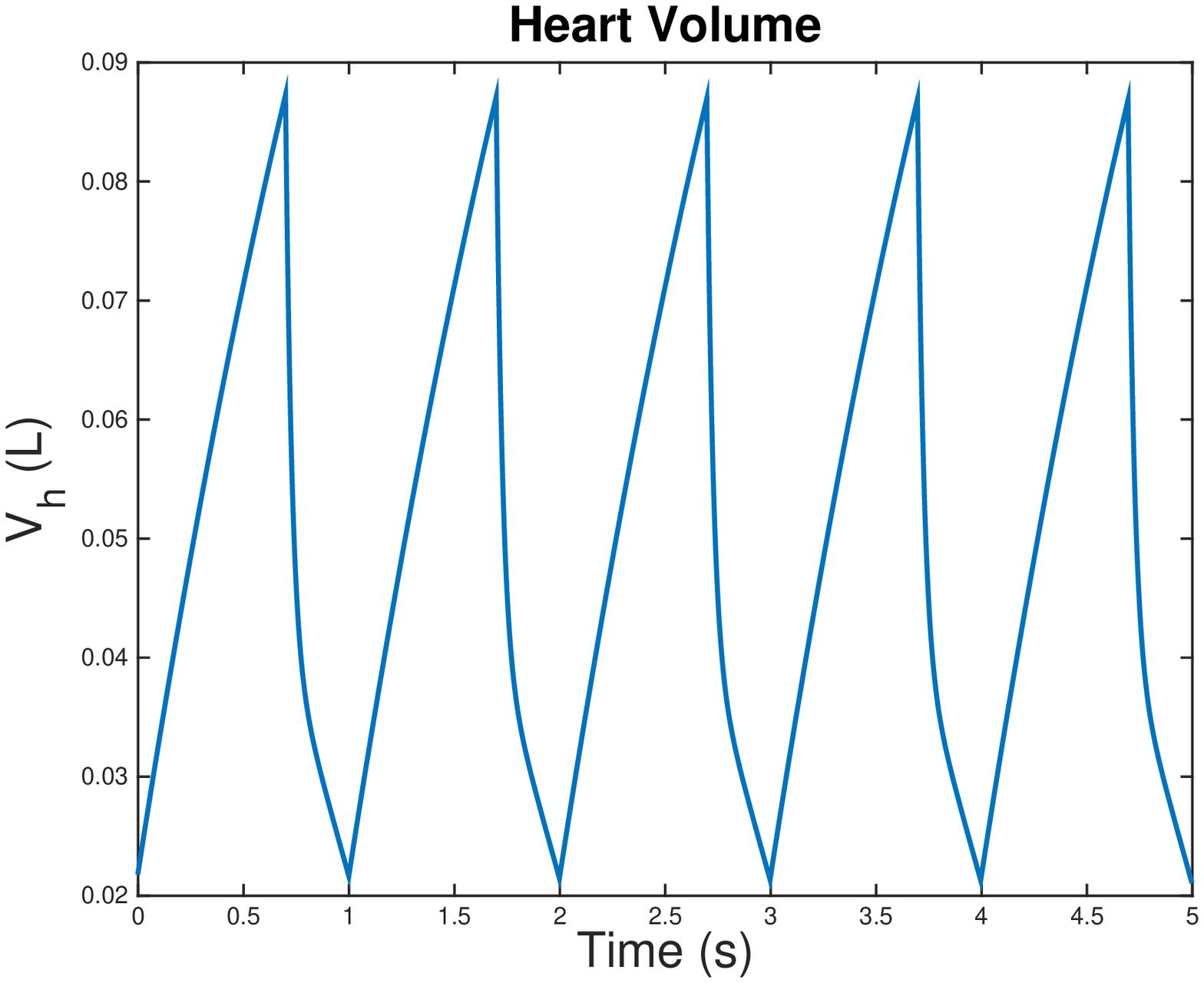}
\hspace{1cm}\includegraphics[height=5.4cm]{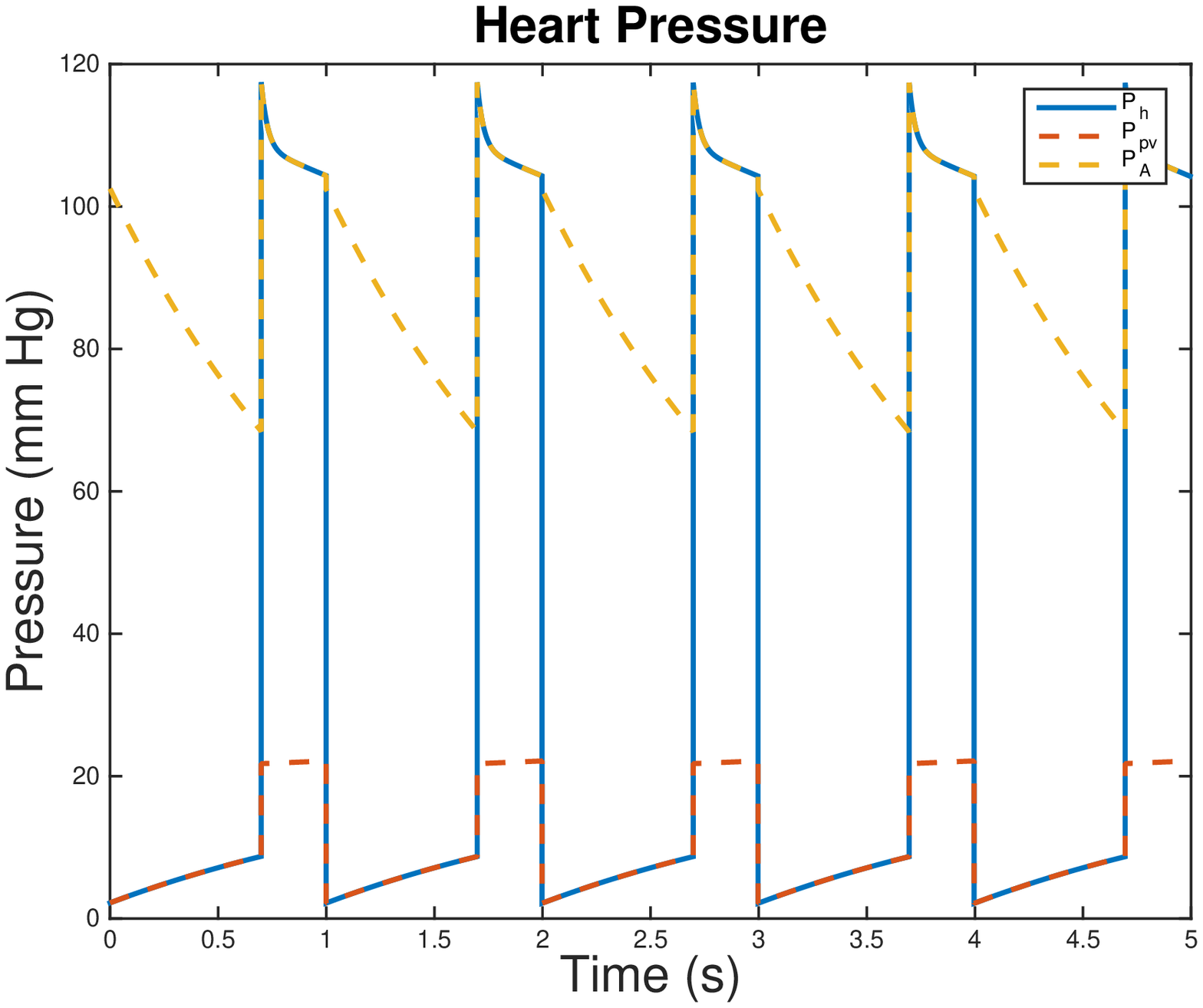}}
\caption{Heart volume (left) and heart pressure (right) as a function of time under normal Fontan conditions.}
\label{fig:fig3}
\end{figure}

Simulations can be done by sequentially integrating the systolic and diastolic equations, and repeating. With model parameter values taken from the table, the simulations of the model equations exhibited realistic values. As shown in Fig.~\ref{fig:fig3}, the stroke volume (the amount of blood pumped out of the heart in one heartbeat) was found to be approximately 70 mL, which is consistent with typical values. The aortic pressure $P_a$ varied between approximately $70$ mm Hg (diastolic) and $120$ mm Hg (systolic) within a typical period of $1$ second, as shown in Fig.~\ref{fig:fig3}.  These values are within the physiological range and the trend is consistent with physiological expectations.  The pulmonary venous pressure $P_{pv}$, which is a surrogate for the atrial pressure, varies between approximately $2.5$ mm Hg and $22.1$ mm Hg.  These values are also within the physiological range.  The trends are reasonably consistent with expectations given our model assumptions.

\begin{figure}
\centerline{ \includegraphics[height=5.4cm]{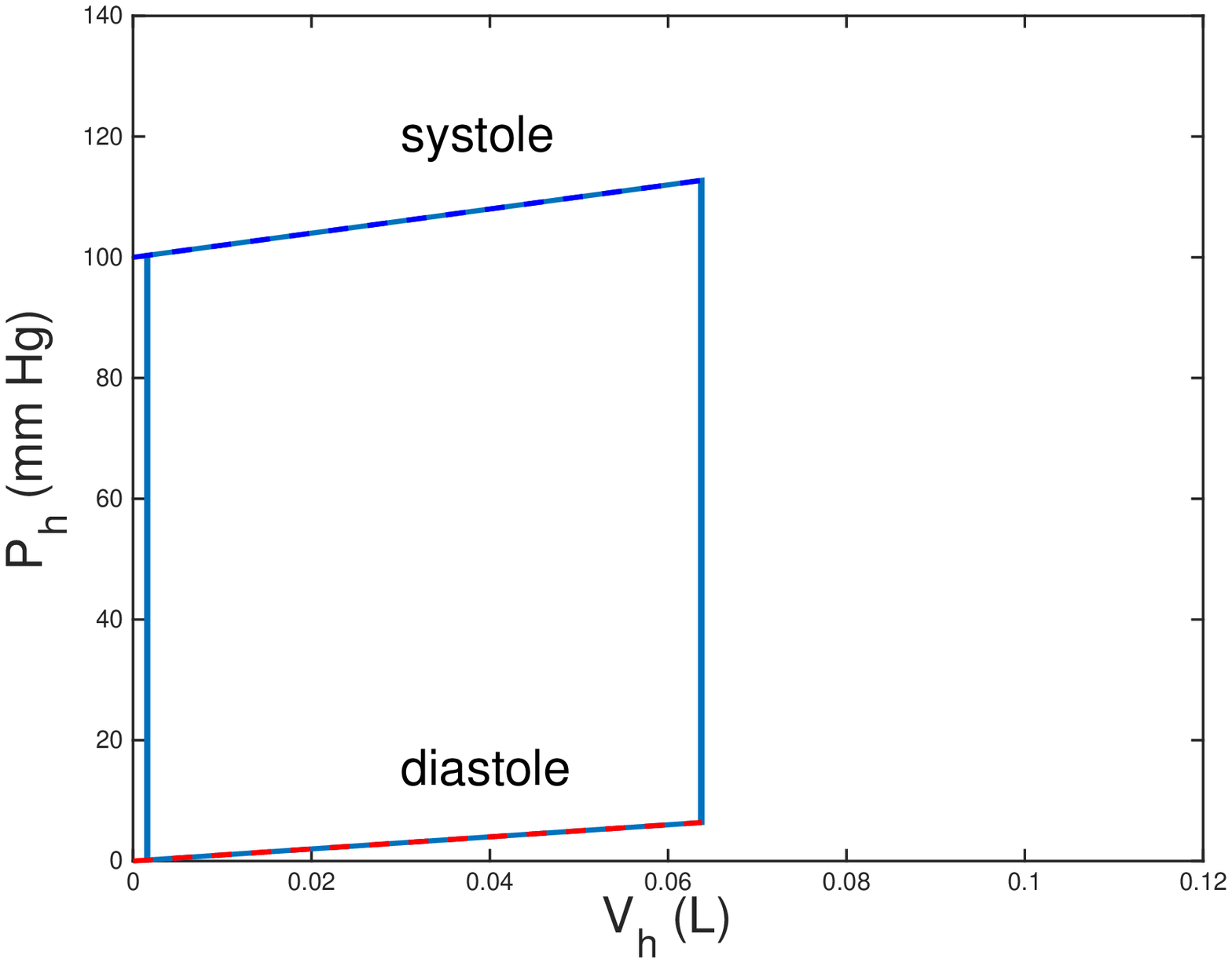}\hspace{1cm}\includegraphics[height=5.4cm]{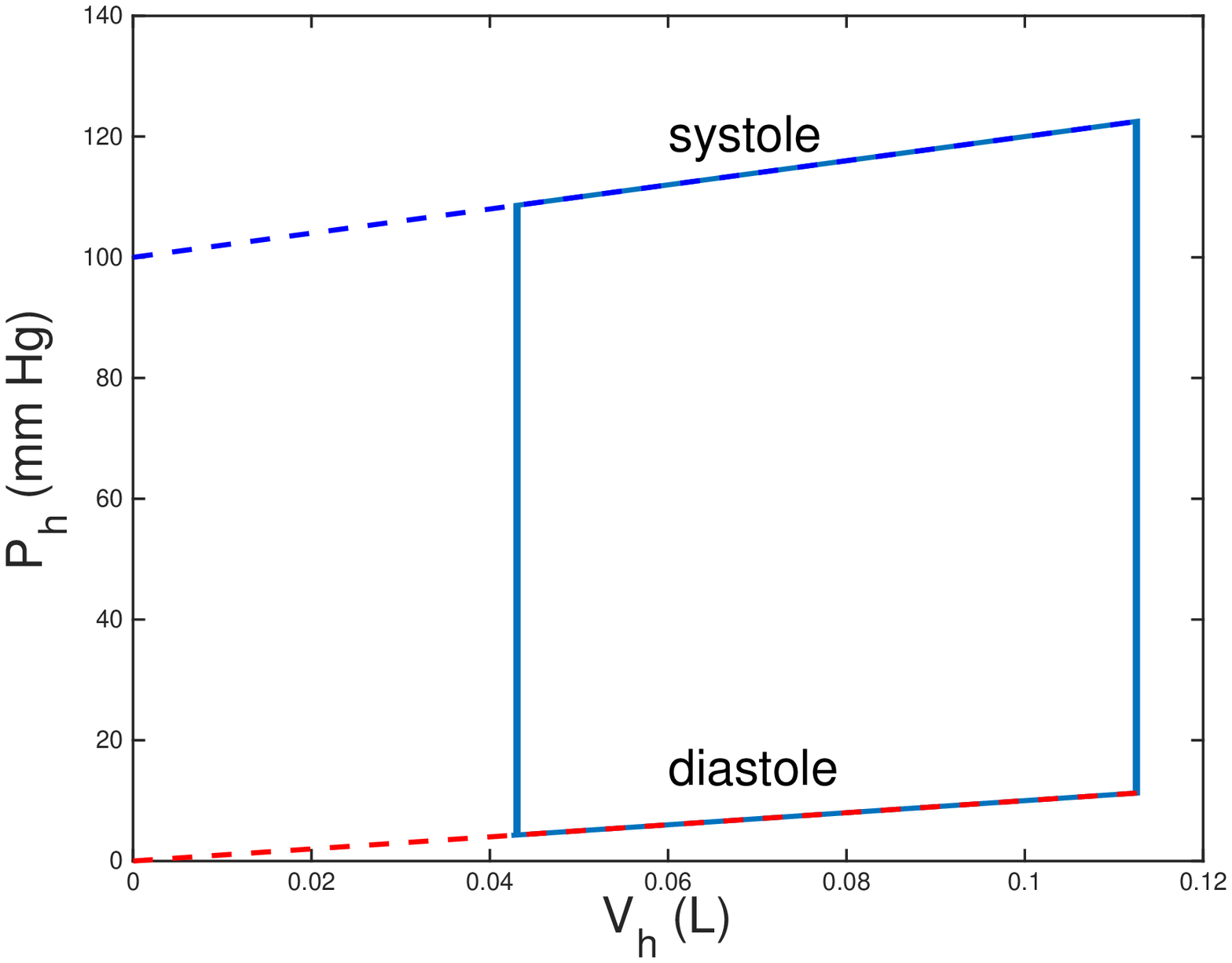} }
\caption{The cardiac pressure-volume cycle from high pulmonary resistance (left) and high heart rate (right).}
\label{fig:figa}
\end{figure}

\begin{figure}
\centerline{\includegraphics[height= 4.5cm]{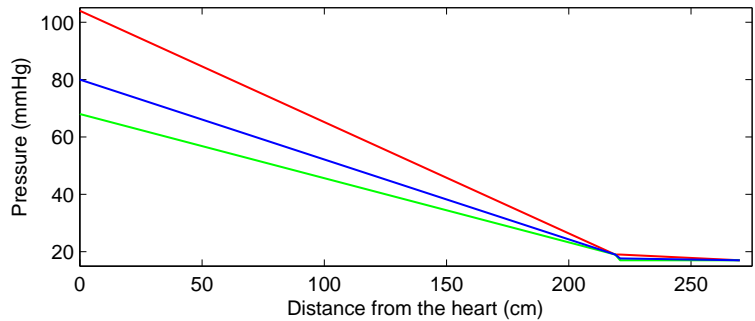}}
\caption{ Blood pressure values (systolic (red), diastolic (green), mean (blue)) for a healthy Fontan patient as a function of distance from the heart.}
\label{fig:figPVR}
\end{figure}
Pulmonary vascular resistance $R_p$ is known to increase in Fontan failure and an increase in this resistance is known to lead to a decrease in cardiac output. This model can be used to demonstrate the impact of pulmonary vascular resistance on cardiac output. Figure~\ref{fig:fig2} shows the change in average cardiac output as a function on pulmonary vascular resistance for two different heart rates.  As expected, the average cardiac output decreases with increasing resistance.  What is interesting is that there is a change in slope of the curves at an inflection point corresponding to $R_p\approx3.60$ mm Hg/min/L for a heart rate of $60$ beats/min and $R_p\approx3.65$ mm Hg/min/L for a heart rate of $120$ beats/min.  After this inflection point, the cardiac output decreases more quickly for a given change in resistance than before this inflection point, suggesting that something changes at this point with regards to Fontan failure and this change is dependent on heart rate.

Pulmonary resistance also has an impact on the cardiac pressure-volume curve. As shown in the left panel of Fig.~\ref{fig:figa}, the cardiac pressure-volume curve shifts to the left (i.e. decreased cardiac volumes) for the case of high pulmonary resistance. What this means is that the basal volume of blood in the heart has decreased as a result of this increase in resistance.  For the present case, the basal volume of the heart has decreased to almost zero, implying that a further increase in resistance would result in an insufficient amount of blood returning to the heart, which consequently would reduce the cardiac output.  Conversely, we see the opposite effect in the right panel Fig.~\ref{fig:figa}, which shows the cardiac pressure-volume curve shifted to the right (i.e. increased cardiac volumes) for the case of high heart rate. By increasing the heart rate for a fixed stroke volume, the cardiac output would increase, resulting in an increase in the amount of blood being pumped to the body and returning to the heart.

Figure~\ref{fig:figPVR} shows the pressure drop as a function of distance from the heart for healthy and failing Fontan patients based on clinically measured pressure catheter data.  This figure illustrates that the majority of the pressure drop occurs near the heart in the systemic arteries and that furthest away from the heart, in the total cavopulmonary connection and the pulmonary circulation, the pressures are low and nearly constant.

\section{Spatially inhomogeneous PDE models of blood pressure distribution}
In this section, our ODE approach to modelling blood flow in the Fontan circulation is extended to a PDE model.  The PDE model has the advantage of allowing for spatial variation of model parameters such as compliance and resistance. This will allow for more personalization of the model to an individual patient, which should improve the accuracy and predictive capabilities of the model.  Furthermore, with piecewise constant values of the model parameters, the PDE model should show a similar behavior to the ODE model, giving us some assurances as to the fidelity of the PDE approach.

To model the circulatory system as a continuous flow network in a resistive compliance vessel, we assume that blood flow is a Stokes flow, i.e. the Reynolds number is sufficiently small to allow us to neglect inertial effects. Consequently, the flux in a cylindrical tube is a Poiseuille flow for which
\begin{equation}
Q = -\tfrac{P_x}{8 \pi \mu} A^2,
\end{equation}
where $P$ is the local pressure, $A$ is the cross-sectional area, and $\mu$ is the fluid viscosity. Now we assume that
a vessel is a linear compliance vessel with $A = A_0 + C\,P$, where $C$ is the compliance. This leads to a flux relationship
for a single vessel
\begin{equation}
Q = -\tfrac{P_x}{8 \pi \mu} (A_0 + C\,P)^2.
\end{equation}
If we have a total number of $N$ parallel vessels all with cross-section area $A$, the flux is
\begin{equation}
Q = -\tfrac{P_x}{8 \pi \mu} \, N \, (A_0 + C\,P)^2 \equiv -q(x,P)\,P_x.
\end{equation}
Notice that in the limit $C \rightarrow 0$, this reduces to Ohm's Law (as it must)
\begin{equation}
Q = -\tfrac{P_x}{R},
\end{equation}
where $R = \tfrac{8 \pi \mu}{N\,A_0^2}$ is the resistance per unit length.
When combined with the conservation law (the total volume of circulating blood is conserved)
\begin{equation}
\label{conserv}
\tfrac{\partial A}{ \partial t} + \tfrac{\partial Q}{\partial x} = 0,
\end{equation}
this yields
\begin{equation}
C\, \tfrac{\partial P}{\partial t} = \tfrac{\partial}{\partial x} \Bigl( q(x,P) \, \tfrac{\partial P}{\partial x}\Bigr),
\end{equation}
which is a nonlinear parabolic partial differential equation for $P(x,t)$ with spatially variable coefficients. In general $C = C(x) \geq 0$, $A_0 = A_0(x) \geq \inf(A_0) > 0$,
and $N = N(x) \geq 1$.

\subsection{Boundary conditions for normal circulation}
\begin{figure}
\centerline{\includegraphics[height=6cm]{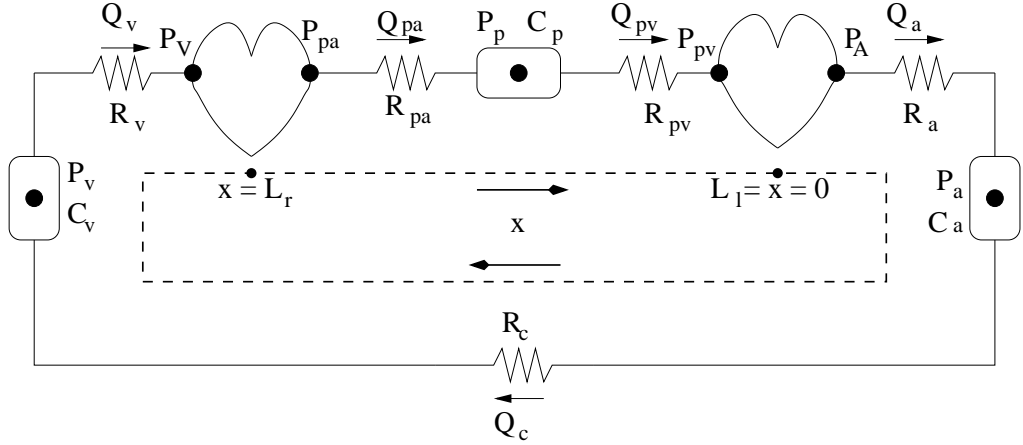}}
\caption{Circuit diagram for the normal circulation.}\label{fig:fig7}
\end{figure}
For the normal circulation (see  Fig.~\ref{fig:fig7}), there are two hearts, the left and the right hearts, each of which satisfy a volume-compliance relationship of the form
\begin{equation}
V = V_0 + C\, P.
\end{equation}
The basal volumes and compliances are different during systole and diastole periods, that is
\begin{equation}\label{eq:left}
V_{lh}(t) = \left\{
           \begin{array}{ll}
             V_{ld}^0 + C_{ld}\,P_{pv}(L_l^{-},t), & \hbox{during diastole} \\
             V_{ls}^0 + C_{ls}\,P_{A}(0^{+},t), & \hbox{during systole}
           \end{array}
         \right.,
\end{equation}
\begin{equation}\label{eq:right}
V_{rh}(t) = \left\{
           \begin{array}{ll}
             V_{rd}^0 + C_{rd}\,P_{V}(L_r^{-},t), & \hbox{during diastole} \\
             V_{rs}^0 + C_{rs}\,P_{pa}(L_r^{+},t), & \hbox{during systole}
           \end{array}
         \right..
\end{equation}
Note that in equations (\ref{eq:left}) and (\ref{eq:right}), the subscripts $l$ and $r$ refer to the left and right hearts, the subscripts $d$ and $s$ refer to diastole and systole. In this two heart model, the 7 pressures are systemic venous pressure $P_v$, vena cava pressure $P_V$, pulmonary arterial pressure $P_{pa}$, pulmonary pressure $P_p$, pulmonary venous pressure $P_{pv}$, aortic pressure $P_A$, and systemic arterial pressure $P_a$.

During systole the input valves (mitral and tricuspid) are closed and output valves (pulmonary and aortic)
are open, while during diastole the opposite is the case. For the normal circulation, we let $0 < x < L_r$
be the systemic circulation, and $L_r < x < L_l$
be the pulmonary circulation, and of course the domain $0 < x < L_l$ is periodic.

For convenience, we use the following notation, during systole: systemic pressure is $P_s^1 := P_{A}$, pulmonary pressure is $P_p^1 := P_{pa}$
and  during diastole: systemic pressure is  $P_s^2 := P_{pv}$, pulmonary pressure is $P_p^2 := P_{V}$. Systemic and pulmonary pressures are everywhere continuous functions
except two points $x=L_l=0$ and $x=L_r$ where discontinuity jumps correspond to the left and right heart pressure jumps during the switch between systole and diastole phases.

\subsection{Initial-boundary value problem (systolic regime: $0 < t < t_1$).}
Partial differential equation:
\begin{equation}
\label{pde}
C(x)\,\tfrac{\partial P}{\partial t} = q(x,P) \, \tfrac{\partial P}{\partial x}.
\end{equation}
Assume that at the initial time $t = 0$ the pressure is $P(x,0)=P_0(x)$. Systemic circulation model $P_s^1$ is given by the partial differential equation (\ref{pde}) and conditions:
\begin{equation}
\left\{
  \begin{array}{ll}
    \hbox{Problem} \,\,\, P_s^1: & 0 < t < t_1, \,\,\, 0 < x < L_r, \\
    \hbox{Initial data}: & P^1_{s\,0}(x) = P_0(x), \\
    \hbox{Boundary conditions for} \,\,\, P_s^1: & B_s^1(t) : = P_s^1(0,t), \\
    \hbox{Dynamic Flux BC at} \,\,\, x = 0: &  C_{ls} \, \tfrac{dB_s^1}{dt} = q(0,B_s^1) \, \tfrac{\partial P_s^1}{\partial x} (0,t), \\
    \hbox{Neumann BC at} \,\,\, x = L_r: & P^1_{s\,x}(L_r,t) = 0. \\
  \end{array}
\right.
\end{equation}

Pulmonary circulation model $P_p^1$ is given by the partial differential equation (\ref{pde}) and conditions:
\begin{equation}
\left\{
  \begin{array}{ll}
    \hbox{Problem} \,\,\, P_p^1: & 0 < t < t_1, \,\,\, L_r < x < L_l, \\
    \hbox{Initial data}: & P^1_{p\,0}(x) = P_0(x), \\
    \hbox{Boundary conditions for} \,\,\, P_p^1: & B_p^1(t) : = P_p^1(L_r,t), \\
    \hbox{Dynamic Flux BC at} \,\,\, x = L_r: &  C_{rs} \, \tfrac{dB_p^1}{dt} = q(L_r,B_p^1) \tfrac{\partial P_p^1}{\partial x} (L_r,t), \\
    \hbox{Neumann BC at} \,\,\, x = L_l: & P^1_{p\,x}(L_l,t) = 0. \\
  \end{array}
\right.
\end{equation}

 These boundary conditions are defined in such a way as to ensure conservation of the total volume
 $$V_T= V_{lh}(t) + V_{rh}(t) + \int_0^{L_l} A(x,t)\,dx$$
 during systole. Indeed, from  (\ref{conserv}), (\ref{eq:left}), (\ref{eq:right}) and the dynamical boundary conditions it follows directly that $\tfrac{d}{dt} V_T = 0$.

\subsection{Initial-boundary value problem  (switch from systole to diastole at: $t = t_1$).}

Initial data for the systemic circulation model $P_s^2$  (and for $B_s^2(t) : = P_s^2(L_r,t)$) at the beginning of the diastole:
\begin{equation}
\left\{
  \begin{array}{ll}
    V_{ls}^0 + C_{ls}\,\, B_s^1(t_1) = V_{ld}^0 + C_{ld} \,\, B_s^2(t_1) & \hbox{conservation of volume,} \\
    P_s^2(x, t_1) = P_s^1(x, t_1), \,\, \hbox{for} \,\,\, x \in (0, L_r - \epsilon) & \hbox{initial data,}\\
    P_s^2(L_r, t_1) = B_s^2(t_1), \,\, \hbox{and} \,\,  P_s^2(x, t_1) \,\, \hbox{is smooth for} \,\, x \in (L_r - \epsilon, L_r) & \hbox{interpolation.} \\
  \end{array}
\right.
\end{equation}
Initial data for the pulmonary circulation model $P_p^2$  (and for $B_p^2(t) : = P_p^2(L_l,t)$) at the beginning of the diastole:
\begin{equation}
\left\{
  \begin{array}{ll}
    V_{rs}^0 + C_{rs}\,\, B_p^1(t_1) = V_{rd}^0 + C_{rd} \,\, B_p^2(t_1) & \hbox{conservation of volume,} \\
    P_p^2(x, t_1) = P_p^1(x, t_1), \,\, \hbox{for} \,\,\, x \in (L_r, L_l - \epsilon) & \hbox{initial data,}\\
    P_p^2(L_l, t_1) = B_p^2(t_1), \,\, \hbox{and} \,\,  P_p^2(x, t_1) \,\, \hbox{is smooth for} \,\, x \in (L_r, L_l - \epsilon) & \hbox{interpolation.} \\
  \end{array}
\right.
\end{equation}

\subsection{Initial-boundary value problem  (diastolic regime: $t_1 < t < t_2$).}

Systemic circulation model $P_s^2$ is given by the partial differential equation (\ref{pde}) and conditions:
\begin{equation}
\left\{
  \begin{array}{ll}
    \hbox{Problem} \,\,\, P_s^2: & t_1 < t < t_2, \,\,\, 0 < x < L_r, \\
    \hbox{Initial data}: & P^2_s (x, t_1), \\
    \hbox{Boundary conditions for} \,\,\, P_s^2: & B_s^2(t) : = P_s^2(L_r,t), \\
    \hbox{Neumann BC at} \,\,\, x = 0: & P^2_{s\,x}0,t) = 0, \\
    \hbox{Dynamic Flux BC at} \,\,\, x = L_r: &  C_{ld} \, \tfrac{dB_s^2}{dt} = q(L_r, B_s^2) \, \tfrac{\partial P_s^2}{\partial x} (L_r,t). \\
  \end{array}
\right.
\end{equation}
Pulmonary circulation model $P_p^2$ is given by the partial differential equation (\ref{pde}) and conditions:
\begin{equation}
\left\{
  \begin{array}{ll}
    \hbox{Problem} \,\,\, P_p^2: & t_1 < t < t_2, \,\,\, L_r < x < L_l, \\
    \hbox{Initial data}: & P^1_p (x, t_1), \\
    \hbox{Boundary conditions for} \,\,\, P_p^2: & B_p^2(t) : = P_p^2(L_l,t), \\
    \hbox{Neumann BC at} \,\,\, x = L_r: & P^2_{p\,x}(L_r,t) = 0, \\
    \hbox{Dynamic Flux BC at} \,\,\, x = L_l: &  C_{rd} \, \tfrac{dB_p^2}{dt} = q(L_l,B_p^2)\, \tfrac{\partial P_p^2}{\partial x} (L_l,t). \\
  \end{array}
\right.
\end{equation}

Again, it follows directly from  (\ref{conserv}), (\ref{eq:left}), (\ref{eq:right}) and the dynamical boundary conditions that during diastole $\tfrac{d}{dt} V_T = 0$.

\subsection{Initial-boundary value problem  (switch from diastole to systole  regime: $t = t_2$).}

Initial data for the systemic circulation model $P_s^1$  (and for $B_s^1(t) : = P_s^1(0,t)$) at the beginning of the systole:
\begin{equation}
\left\{
  \begin{array}{ll}
    V_{ls}^0 + C_{ls}\,\, B_s^1(t_2) = V_{ld}^0 + C_{ld} \,\, B_s^2(t_2) & \hbox{conservation of volume,} \\
    P_s^1(x, t_2) = P_s^2(x, t_2), \,\, \hbox{for} \,\,\, x \in (\epsilon, L_r) & \hbox{initial data,}\\
    P_s^1(0, t_2) = B_s^1(t_2), \,\, \hbox{and} \,\,  P_s^1(x, t_2) \,\, \hbox{is smooth for} \,\, x \in (0, \epsilon) & \hbox{interpolation.} \\
  \end{array}
\right.
\end{equation}
Initial data for the pulmonary circulation model $P_p^1$  (and for $B_p^1(t) : = P_p^1(L_r,t)$) at the beginning of the systole:
\begin{equation}
\left\{
  \begin{array}{ll}
    V_{rs}^0 + C_{rs}\,\, B_p^1(t_2) = V_{rd}^0 + C_{rd} \,\, B_p^2(t_2) & \hbox{conservation of volume,} \\
    P_p^1(x, t_2) = P_p^2(x, t_2), \,\, \hbox{for} \,\,\, x \in (L_r + \epsilon, L_l) & \hbox{initial data,}\\
    P_p^1(L_r, t_2) = B_p^1(t_2), \,\, \hbox{and} \,\,  P_p^1(x, t_2) \,\, \hbox{is smooth for} \,\, x \in (L_r, L_r +\epsilon) & \hbox{interpolation.} \\
  \end{array}
\right.
\end{equation}
The problem is periodic in time (i.e systolic and diastolic regimes are repeated).

\subsection{Boundary conditions for Fontan circulation}

The Fontan blood flow circulation has only one heart, so the model has only a single loop $0 < x < L_l$.
The basal volumes and compliances of the heart are also different during systole and diastole periods, that is
\begin{equation}
V_{h} = \left\{
           \begin{array}{ll}
             V_{d}^0 + C_{d}\,P_{pv}, & \hbox{during diastole} \\
             V_{s}^0 + C_{s}\,P_{A}, & \hbox{during systole}
           \end{array}
         \right..
\end{equation}

For convenience we use the following notation, during systole: blood pressure is $P^1 := P_{A}$
and  during diastole blood pressure is  $P^2 := P_{pv}$.

\subsection{Initial-boundary value problem (systolic regime: $0 < t < t_1$).}

Assume that at the initial time $t = 0$ the pressure is $P(x,0) =P_0(x)$.

One heart circulation model $P^1$ is given by the partial differential equation (\ref{pde}) and conditions:
\begin{equation}
\left\{
  \begin{array}{ll}
    \hbox{Problem} \,\,\, P^1: & 0 < t < t_1, \,\,\, 0 < x < L_l, \\
    \hbox{Initial data}: & P^1_{0}(x) = P_0(x), \\
    \hbox{Boundary conditions for} \,\,\, P^1: & B^1(t) : = P^1(0,t), \\
    \hbox{Dynamic Flux BC at} \,\,\, x = 0: &  C_{ls} \, \tfrac{dB^1}{dt} = q(0,B^1) \, \tfrac{\partial P^1}{\partial x} (0,t), \\
    \hbox{Neumann BC at} \,\,\, x = L_l: & P^1_{x}(L_l,t) = 0. \\
  \end{array}
\right.
\end{equation}

\subsection{Initial-boundary value problem  (switch from systole to diastole at: $t = t_1$).}

Initial data for the systemic circulation model $P^2$  (and for $B^2(t) : = P^2(L_l,t)$) at the beginning of the diastole:
\begin{equation}
\left\{
  \begin{array}{ll}
    V_{s}^0 + C_{s}\,\, B^1(t_1) = V_{d}^0 + C_{d} \,\, B^2(t_1) & \hbox{conservation of volume,} \\
    P^2(x, t_1) = P^1(x, t_1), \,\,  \hbox{for} \,\,\, x \in (0, L_l - \epsilon) & \hbox{initial data,}\\
    P^2(L_l, t_1) = B^2(t_1), \,\, \hbox{and} \,\,  P^2(x, t_1) \,\, \hbox{is smooth for} \,\, x \in (L_l - \epsilon, L_l) & \hbox{interpolation.} \\
  \end{array}
\right.
\end{equation}

\subsection{Initial-boundary value problem  (diastolic regime: $t_1 < t < t_2$).}

Systemic circulation model $P^2$ is given by the partial differential equation (\ref{pde}) and conditions:
\begin{equation}
\left\{
  \begin{array}{ll}
    \hbox{Problem} \,\,\, P^2: & t_1 < t < t_2, \,\,\, 0 < x < L_l, \\
    \hbox{Initial data}: & P^2(x, t_1), \\
    \hbox{Boundary conditions for} \,\,\, P^2: & B^2(t) : = P^2(L_l,t), \\
    \hbox{Neumann BC at} \,\,\, x = 0: & P^2_{x}0,t) = 0, \\
    \hbox{Dynamic Flux BC at} \,\,\, x = L_l: &  C_{ld} \, \tfrac{dB^2}{dt} =  q(L_l,B^2) \tfrac{\partial P^2}{\partial x} (L_l,t).\\
  \end{array}
\right.
\end{equation}

\subsection{Initial-boundary value problem  (switch from diastole to systole  regime: $t = t_2$).}

Initial data for the systemic circulation model $P^1$  (and for $B^1(t) : = P^1(0,t)$) at the beginning of the systole:
\begin{equation}
\left\{
  \begin{array}{ll}
    V_{s}^0 + C_{s}\,\, B^1(t_2) = V_{d}^0 + C_{d} \,\, B^2(t_2), & \hbox{conservation of volume,} \\
    P^1(x, t_2) = P^2(x, t_2), \,\, \hbox{for} \,\,\, x \in (\epsilon, L_r), & \hbox{initial data,}\\
    P^1(0, t_2) = B^1(t_2), \,\, \hbox{and} \,\,  P^1(x, t_2) \,\, \hbox{is smooth for} \,\, x \in (0, \epsilon) & \hbox{interpolation.} \\
  \end{array}
\right.
\end{equation}

The problem is periodic in time (i.e systolic and diastolic regimes are repeated). It is a direct computation to verify that with these boundary conditions $\tfrac{d}{dt} V_T = 0$.

\section{Well-posedness analysis of spatially inhomogeneous PDE model}

\subsection{Existence and uniqueness of the nonnegative solution in normal case.}
In this section we find restrictions on the parameters of the PDE model
with dynamical boundary conditions for the case of the normal blood circulation
for which a nonnegative solution exists and is unique. First, we show that the pressure $P(x,t)$ stays bounded for any time $t$, then we obtain uniform in time bounds for the time derivative and for the gradient of $P(x,t)$ for some range of parameter values.
Second, we employ Gr\"{o}nwall lemma to prove uniqueness.

We start by introducing the following notations
$$
Q_{t_1}^{t_2} = \Omega_r \times (t_1,t_2),\  Q_{t_2} = \Omega_r \times (0,t_2),
\ \Omega_r := (0,L_r), \  \Omega_r^{\epsilon} := ( L_r - \epsilon, L_r).
$$
Consider the following equation
\begin{equation}\label{a-1}
C(x) P_t = (q(x,P) \, P_x)_x  
\end{equation}
with initial and boundary conditions

\begin{equation}\label{a-4}
P_1(x,0) = P_0(x)  \geqslant 0,
\end{equation}
\begin{equation}\label{a-2-1}
C_{ls}  P_{1,t} =  q(x, \, P_1) P_{1,x}  \text{ at } x = 0,\
P_{1,x} = 0 \text{ at } x = L_r
\end{equation}
for all $t \in (0,t_1)$ (during systole);

\begin{equation}\label{a-4-2}
P_2(x,t_1) = \left \{
\begin{gathered}
P_1(x,t_1) \ \forall\,x \in [0,L_r -\epsilon],  \\
f(x,t_1) \ \forall\,x \in (L_r -\epsilon,L_r],
\end{gathered}  \right.
\end{equation}
where $f(x,t_1)$ is such that
\begin{equation}\label{f}
\int \limits_{\Omega_r^{\epsilon}} { C(x)
f(x,t_1) dx} - C_{ld} f(L_r,t_1) =
\int \limits_{\Omega_r^{\epsilon}} { C(x) P_1 (x,t_1) dx}
+ C_{ls} P_1(0,t_1),
\end{equation}
\begin{equation}\label{a-2-2}
P_{2,x} = 0 \,\, \text{ at } x = 0, \,\,  C_{ld} P_{2,t} = q(x,P_2) \, P_{2,x} \,\, \text{ at } \,\, x = L_r
\end{equation}
for all $t \in (t_1,t_2)$ (during diastole);

\begin{equation}\label{a-4-2-00}
P_1(x,t_2) = \left \{
\begin{gathered}
\tilde{f}(x,t_2) \ \forall\,x \in [0, \epsilon],  \\
P_2(x,t_2) \ \forall\,x \in ( \epsilon,L_r],
\end{gathered}  \right.
\end{equation}
where $\tilde{f}(x,t_2)$ is such that
\begin{equation}\label{f-00}
\int \limits_{0}^{\epsilon} { C(x)
\tilde{f}(x,t_2) dx} + C_{ls} \tilde{f}(0,t_2) =
\int \limits_{0}^{\epsilon} { C(x) P_2 (x,t_2) dx}
- C_{ld} P_2(L_r,t_2).
\end{equation}

Let
\begin{equation}\label{sol}
P(x,t) = \left \{
\begin{gathered}
P_1(x,t) \ \forall\,t \in (0,t_1)  \\
P_2(x,t) \ \forall\,t \in [t_1,t_2)
\end{gathered}  \right.
\end{equation}
be a solution to problem (\ref{a-1})--(\ref{a-2-2}) for all $t \in (0,t_2)$ such that
\begin{equation}\label{period}
P(x,t) = P(x,t+t_2).
\end{equation}

First, we define a weak solution for our problem.
\begin{definition}\label{def}
A non-negative function $P(x,t)$ is said to be a periodic solution of the problem (\ref{a-1})--(\ref{a-2-2}),
i.\,e. $P(x,t) = P(x,t+t_2)$, if
$$
P \in C_{x,t}^{\tfrac{1}{2},\tfrac{1}{4}}(\bar{Q}_{t_2}) \cap  L^{\infty}(0,t_2;H^1(\Omega_r)) ,\ C^{\tfrac{1}{2}}(x)P_t \in L^2(Q_{t_2}),
$$
$$
q(x,P) \, P_x  \in L^2 (0,t_2;H^1(\Omega_r)),
$$
and $P(x,t)$ satisfies equation (\ref{a-1})  in
the sense that
$$
\iint \limits_{Q_{t_2}} {C(x) P_t \phi  \,dx dt} -  \iint \limits_{Q_{t_2}}
{( q(x,P) \, P_x )_x \phi  \,dx dt} = 0
$$
for any $\phi \in L^2(Q_{t_2})$ and  $\phi(x,0) = \phi(x,t_2)$.
\end{definition}

Our main result establishes parameter ranges for which non-negative solutions exist, as follows.
\begin{theorem} \label{ex}
If
$$
C(x) \in C(\bar{\Omega}_r) :  \inf C(x) > 0, \ \|C(x)\|_1   < C_{ld}; \ C_{ls} \geqslant 0,  \ C_{ld} > 0;
$$
$$
\inf N(x) > 0, \,\,  \inf A_0(x) > 0,\  f(x,t_1) \in C(\bar{\Omega}_r^\epsilon),
$$
initial data $P_0(x) \in H^1(\Omega_r)$ are non-negative, and
$$
\int \limits_{\Omega_r} {  q(x,P_0) \, P_{0,x}^2 dx} \leqslant
\tfrac{1}{64 \pi \mu^2 |\Omega_r|}\tfrac{(C_{ld} - \| C(x)\|_1)}{C^2_{ld}} \mathop {\inf}\limits_{\Omega_r}| \tfrac{N(x) A_0^4(x)}{ C^2(x)}|,
$$
$$
M:= \int \limits_{\Omega_r} { C(x) P_0 (x) dx} + C_{ls}P_0 (0)   > 0 ,
$$
then the problem (\ref{a-1})--(\ref{a-2-2}) has a unique non-negative solution
in the sense of the Definition~\ref{def}.
\end{theorem}

\subsection{Proof of Theorem~\ref{ex}}\label{proof1}
Note that the equation (\ref{a-1}) becomes degenerate if $P = - \inf \frac{A_0(x)}{C(x)}$. Hence, we start by constructing a sequence of positive approximations of initial data $P_{0n}>0$. We can choose for example $P_{0n}(x) = P_0(x)+ \frac{1}{n}$ (if $n \to \infty$ then $P_{0n}(x) \to P_0(x)$). These approximations allow us to apply the theoretical background developed for uniformly parabolic equations to our problem. By taking $n \to \infty$, as a limit, we obtain a weak solution $P(x,t)$.  We omit some technical details and derive only a priori estimates which imply the existence of this weak solution. We need to specify conditions on the smoothing functions $f(x,t)$ and $\tilde{f}(x,t)$ such that total volume is conserved on the whole time interval $[0,t_2]$.

\textbf{\emph{Volume conservation}:} Integrating  (\ref{a-1}) on $Q_{t}$, due to
(\ref{a-2-1}) and (\ref{a-4}), we arrive at
\begin{multline} \label{a-00}
 \int \limits_{\Omega_r} { C(x) P_1 (x,t) dx} +
C_{ls} P_1 (0,t) = \\
\int \limits_{\Omega_r} { C(x) P_0 (x) dx} + C_{ls}P_0 (0) =: M > 0 \ \forall\,t\in [0,t_1).
\end{multline}
Integrating  (\ref{a-1}) on $Q_{t_1}^{t}$, due to
(\ref{a-2-2}) and (\ref{a-4-2}), we have
\begin{multline} \label{a-00-0}
 \int \limits_{\Omega_r} { C(x) P_2 (x,t) dx} +
C_{ld} P_{2} (L_r,t_1) = \\
\int \limits_{\Omega_r} { C(x) P_2(x,t_1) dx} + C_{ld}P_2 (L_r,t) \ \forall\,t\in [ t_1,t_2).
\end{multline}
By (\ref{a-00}) and (\ref{a-00-0}), we obtain
\begin{multline} \label{a-00-1}
 \int \limits_{\Omega_r} { C(x) P_2 (x,t) dx} + \int \limits_{\Omega_r^{\epsilon}} { C(x)
(P_1 (x,t_1) - f(x,t_1) ) dx} +
C_{ld}(f (L_r,t_1) - P_2(L_r,t))  = \\
\int \limits_{\Omega_r} { C(x) P_0(x) dx} +
 C_{ls}(P_0(0) - P_{1} (0,t_1)) .
\end{multline}
Due to (\ref{f})
, from (\ref{a-00-1}) we get
\begin{equation} \label{a-00-2}
 \int \limits_{\Omega_r} { C(x) P_2 (x,t) dx}  -
 C_{ld} P_2(L_r,t) = M \ \forall\,t\in [ t_1,t_2).
\end{equation}
Moreover, by (\ref{a-4-2-00}) and (\ref{a-00-2}) we find that
\begin{multline*}
 \int \limits_{\Omega_r} { C(x) P_1 (x,t_2) dx} +  C_{ls} P_1(0,t_2) = \\
 \int \limits_{0}^{\epsilon} { C(x) \tilde{f}(x,t_2) dx} + \int \limits_{\epsilon}^{L_r} { C(x) P_2(x,t_2)  dx}
 +  C_{ls} \tilde{f}(0,t_2) =\\
\int \limits_{0}^{\epsilon} { C(x)( \tilde{f}(x,t_2) - P_2(x,t_2)) dx} +
C_{ld}P_2(L_r,t_2) +
C_{ls} \tilde{f} (0,t_2)  + M  ,
\end{multline*}
whence, due to (\ref{f-00}), we have
\begin{equation} \label{a-00-33}
 \int \limits_{\Omega_r} { C(x) P_1 (x,t_2) dx} +  C_{ls} P_1(0,t_2) = M.
\end{equation}
Consequently, total volume of the left heart and systolic circulation is identical at
$t = 0$ and $t = t_2$.

Below we show how to prove, using Moser's method \cite{moser1966rapidly}, that the blood pressure $P(x,t)$
stays bounded on the whole time interval $[0,t_2]$. We start by showing that
$P(x,t)$ is bounded in $L^{\infty}(0,t_2; L^2(\Omega_r))$ then we show that for any $\alpha > 0$ the solution
$P(x,t)$ is bounded in $L^{\infty}(0,t_2; L^{\alpha +2}(\Omega_r) )$ and after that we take the limit
$\alpha \rightarrow \infty$.

\noindent \textbf{\emph{Boundedness}:} Multiplying (\ref{a-1}) by $P(x,t)$ and integrating along $\Omega_r$, due to
(\ref{a-2-1}), we have
\begin{equation}\label{b-1}
\tfrac{1}{2}\tfrac{d}{dt} \int \limits_{\Omega_r} { C(x) P^2(x,t) dx} +
 \int \limits_{\Omega_r} { q(x,P) \, P_x^2  dx} =
 -\tfrac{1}{2} \tfrac{d}{dt} [C_{ls} P^2(0,t)].
  \end{equation}
Integrating (\ref{b-1}) in time, we get
\begin{multline}\label{b-2}
\tfrac{1}{2}  \int \limits_{\Omega_r} { C(x) P^2(x,t) dx} + \tfrac{1}{2}C_{ls} P^2(0,t)  +
 \iint \limits_{Q_0^{t_1}} {  q(x,P) \, P_x^2  dx dt} = \\
\tfrac{1}{2}  \int \limits_{\Omega_r} { C(x) P_0^2(x) dx} +
 \tfrac{1}{2}C_{ls} P_0^2(0 ) =: K_0 \ \forall\,t \in (0,t_1) .
\end{multline}
On the other hand, multiplying (\ref{a-1}) by $P(x,t)$ and integrating along $\Omega_r$, due to
(\ref{a-2-2}), we have
\begin{equation}\label{b-1-2}
\tfrac{1}{2}\tfrac{d}{dt} \int \limits_{\Omega_r} { C(x) P^2(x,t) dx} +
 \int \limits_{\Omega_r} {  q(x,P) \, P_x^2  dx} =
 \tfrac{1}{2} \tfrac{d}{dt} [C_{ld} P^2(L_r,t)],
  \end{equation}
and integrating (\ref{b-1-2}) in time from $t_1$, we have
\begin{multline}\label{b-2-2}
\tfrac{1}{2}  \int \limits_{\Omega_r} { C(x) P^2(x,t) dx} + \tfrac{1}{2} C_{ld} P^2(L_r,t_1) +
 \iint \limits_{Q_{t_1}^{t}} {  q(x,P) \, P_x^2  dx dt} = \\
\tfrac{1}{2}  \int \limits_{\Omega_r} { C(x) P^2(x,t_1) dx} +
 \tfrac{1}{2} C_{ld} P^2(L_r,t) \ \forall\,t \in ( t_1, t_2).
\end{multline}
By (\ref{a-00-2}) with $M > 0$ we find that
\begin{equation} \label{f-01}
C_{ld} P^2(L_r,t) \leqslant \tfrac{1}{C_{ld}}\Bigl(\int \limits_{\Omega_r} { C(x) P  (x,t) dx}  \Bigr)^2   \leqslant
\tfrac{\|C(x)\|_1}{C_{ld}} \int \limits_{\Omega_r} { C(x) P^2(x,t) dx}  \ \forall\,t \in ( t_1, t_2).
\end{equation}
As a result, from (\ref{b-2-2}), due to (\ref{f-01}) and (\ref{b-2}), we have
\begin{multline}\label{b-2-3}
\tfrac{1}{2}(1- \tfrac{\|C(x)\|_1}{C_{ld}})  \int \limits_{\Omega_r} { C(x) P^2(x,t) dx} + \tfrac{1}{2} C_{ld} P^2(L_r,t_1) +\\
 \iint \limits_{Q_{t_1}^{t}} {  q(x,P) \, P_x^2  dx dt} \leqslant
\tfrac{1}{2}  \int \limits_{\Omega_r} { C(x) P^2(x,t_1) dx} \leqslant \\
K_1 := K_0 + \int \limits_{\Omega_r^{\epsilon}} { C(x)(f^2(x,t_1) - P_1^2(x,t_1)) dx}  \ \forall\,t \in ( t_1, t_2)
\end{multline}
provided
\begin{equation}\label{con-1}
1- \tfrac{\|C(x)\|_1}{C_{ld}} > 0 \Leftrightarrow 0< \tfrac{\|C(x)\|_1}{C_{ld}} < 1.
\end{equation}

Multiplying (\ref{a-1}) by $P^{\alpha +1}(x,t)$ with $\alpha \geqslant 0$ and integrating along $\Omega_r$
and in time, we have
\begin{multline}\label{b-2-a}
  \int \limits_{\Omega_r} { C(x) P^{\alpha +2}(x,t) dx} +  C_{ls} P^{\alpha +2}(0,t)  +
 (\alpha +1)(\alpha +2) \iint \limits_{Q_0^{t_1}} {  q(x,P) \, P^{\alpha }P_x^2  dx dt} = \\
 \int \limits_{\Omega_r} { C(x) P_0^{\alpha +2}(x) dx} +
  C_{ls} P_0^{\alpha +2}(0 ) =: K_{\alpha} \ \forall\,t \in (0,t_1) ,
\end{multline}
\begin{multline}\label{b-2-b}
(1- (\tfrac{\|C(x)\|_1}{C_{ld}})^{\alpha +1})  \int \limits_{\Omega_r} { C(x) P^{\alpha +2}(x,t) dx}
+   C_{ld} P^{\alpha +2}(L_r,t_1) +\\
(\alpha +1)(\alpha +2) \iint \limits_{Q_{t_1}^{t}} {  q(x,P) \, P^{\alpha} P_x^2  dx dt} \leqslant
   \int \limits_{\Omega_r} { C(x) P^{\alpha +2}(x,t_1) dx} \leqslant \\
K_{\alpha +1} := K_{\alpha} + \int \limits_{\Omega_r^\epsilon} { C(x)(f^{\alpha +2}(x,t_1) - P_1^{\alpha +2}(x,t_1)) dx}
\ \forall\,t \in ( t_1, t_2)
\end{multline}
provided (\ref{con-1}). Next by Moser's method \cite{moser1966rapidly}, taking into account that
$$
\mathop{\sup} \limits_{\Omega_r} |P| = \mathop{\lim} \limits_{\gamma \to +\infty}
\Bigl( \tfrac{1}{|\Omega_r|} \int \limits_{\Omega_r} { P^{\gamma}  dx}  \Bigr)^{\tfrac{1}{\gamma}},
$$
due to the periodicity and $\inf C(x) > 0$, from (\ref{b-2-b}) we obtain
\begin{equation}\label{moser}
\mathop {\sup} \limits_{\Omega_r} |P| \leqslant K < +\infty \ \forall\, t > 0
\end{equation}
 provided $ \mathop {\sup} \limits_{\Omega_r^\epsilon} |f(x,t_1)| < + \infty$
 and $ 0< \tfrac{\|C(x)\|_1}{C_{ld}} < 1$.

Now we obtain the main a priori estimates for the gradient $P_x$ in $L^{\infty}(0,t_2; L^2(\Omega_r) )$ and for the
time derivative $P_t$  in $L^2(Q_{t_2})$. Using these estimates we are able to build a weak solution for the problem at hand.

\noindent \textbf{\emph{A priori estimate}:}  Multiplying (\ref{a-1}) by $P_t$ and integrating along $\Omega_r$, we have
\begin{multline}\label{c-1}
\tfrac{1}{2} \tfrac{d}{dt} \int \limits_{\Omega_r} {q(x,P) \, P_x^2  dx} +
\int \limits_{\Omega_r} { C(x) P_t^2 dx} = \\
\int \limits_{\Omega_r} {  \Bigl( \tfrac{N(x) q(x,P)}{8\pi\mu} \Bigr)^{\tfrac{1}{2}} \, P_x^2 C(x) \, P_t dx} +
 P_t q(x,P) \, P_x \biggl |_{0}^{L_r}.
  \end{multline}
whence, due to (\ref{a-2-1}),
\begin{multline}\label{c-2}
\tfrac{1}{2}  \tfrac{d}{dt}\int \limits_{\Omega_r} {  q(x,P) \, P_x^2 dx} +
 \int \limits_{\Omega_r} { C(x) P_t^2 dx} +
C_{ls}( P(0,t))_t^2 =
\int \limits_{\Omega_r} { ( \tfrac{N q}{ 8 \pi \mu} )^{1/2} P_x^2 C(x) P_t dx} \leqslant \\
\mathop {\sup} \limits_{\Omega_r}| q(x,P) \, P_x |
\Bigl(\int \limits_{\Omega_r} { C(x) P^2_t dx}\Bigr)^{\tfrac{1}{2}}
\mathop {\sup} \limits_{\Omega_r}| \tfrac{C(x)N(x)}{q(x,P)^2} |^{\tfrac{1}{2}}
\Bigl(
\int \limits_{\Omega_r} {  q(x,P) \, P_x^2 dx}\Bigr)^{\tfrac{1}{2}}  \leqslant \\
\mathop {\sup} \limits_{\Omega_r}| \tfrac{C^2(x)|\Omega_r| N(x)}{\pi q(x,P)^2} |^{\tfrac{1}{2}}
\Bigl(\int \limits_{\Omega_r} {  C(x) P^2_t dx}\Bigr)
\Bigl(
\int \limits_{\Omega_r} { q(x,P) \, P_x^2 dx}\Bigr)^{\tfrac{1}{2}}
\end{multline}
for all $t \in (0,t_1)$. Let us denote by
$$
y(t):= \tfrac{1}{2}  \int \limits_{\Omega_r} { q(x,P) \, P_x^2 dx},
\quad
a(t) : =  \int \limits_{\Omega_r} { C(x) P_t^2 dx}, \,\,
b(t) := C_{ls}( P(0,t))_t^2.
$$
Then from (\ref{c-2}) and using that $q \geq \tfrac{1}{8 \pi \mu}  N A^2_0$ we find that
\begin{equation}\label{c-3}
y'(t) + b(t) \leqslant a(t)(C_1 y^{\tfrac{1}{2}}(t) - 1),
\end{equation}
where
$$
C_1 =  \mathop {\sup} \limits_{\Omega_r} \left|\tfrac{ 64 \pi \mu^2 |\Omega_r|  C^2(x) } { N(x) A^4_0(x) } \right|^{\tfrac{1}{2}}.
$$
Indeed, taking into account
\begin{multline*}
\tfrac{1}{4 \pi \mu} \mathop {\inf} \limits_{\Omega_r}(N(x)A_0^2(x))\,y(t)
\leqslant \int \limits_{\Omega_r} {  q(x,P)^2 \, P_x^2  dx} \leqslant    \\
\Bigl(\tfrac{|\Omega_r|}{\pi}\Bigr)^{2} \int \limits_{\Omega_r} { ( q(x,P) \, P_x)_x^2  dx} \leqslant
\Bigl( \tfrac{|\Omega_r|} {\pi} \Bigr)^2 \mathop {\sup} \limits_{\Omega_r}(C(x)) a(t),
\end{multline*}
i.\,e.
\begin{equation}\label{as0}
a(t) \geqslant C_2 y(t), \text{ where } C_2 =  \tfrac{\pi}{4 \mu |\Omega_r|^2}
\mathop {\inf} \limits_{\Omega_r}(\tfrac{N(x)A_0^2(x)}{C(x)}) ,
\end{equation}
from (\ref{c-3}) we get
\begin{equation}\label{c-4}
 y'(t) + b(t) \leqslant C_2y(t)(C_1 y^{\tfrac{1}{2}}(t) -1),
\end{equation}
provided $y(0) < \frac{1}{C_1^2} $. So,
$$
y_1(t) \leqslant \tfrac{\beta^2 y_1(0)}{(y_1^{\tfrac{1}{2}}(0) +(\beta - y_1^{\tfrac{1}{2}}(0) ) e^{\tfrac{\alpha \beta}{2}t})^2}
 \ \forall\, t \in [0,t_1) \text{ if } y(0) < \beta^2,
$$
where $\alpha =C_1 C_2$  and $\beta = \frac{1}{C_1}$.

Similar to (\ref{c-2}), for all $t \in (t_1,t_2)$ we deduce that
\begin{multline}\label{c-2-0}
\tfrac{1}{2} \tfrac{d}{dt}\int \limits_{\Omega_r} {  q(x,P) \, P_x^2 dx} +
 \int \limits_{\Omega_r} { C(x) P_t^2  dx} \leqslant \\
 C_{ld} (P(L_r,t))_t^2  +
\mathop {\sup} \limits_{\Omega_r}| \tfrac{C^2(x)N(x)|\Omega_r|}{\pi q(x,P)^2} |^{\tfrac{1}{2}}
\Bigl(\int \limits_{\Omega_r} {  C(x) P^2_t dx}\Bigr)
\Bigl(
\int \limits_{\Omega_r} { q(x,P) \, P_x^2 dx}\Bigr)^{\tfrac{1}{2}}.
\end{multline}
Let us denote by
$$
d(t) :=  C_{ld} (P(L_r,t))_t^2 .
$$
Then from (\ref{c-2-0}) we find that
\begin{equation}\label{c-3-2}
y'(t) \leqslant a(t)(C_1 y^{\tfrac{1}{2}}(t) - 1) + d(t).
\end{equation}
By the volume conservation (\ref{a-00-2}), we arrive at
\begin{equation}\label{c-000}
 d(t) = \tfrac{1}{C_{ld}} \Bigl( \int \limits_{\Omega_r} {  (C(x) P)_t  dx} \Bigr)^2 \leqslant
\tfrac{\|C(x)\|_1}{C_{ld}}   \int \limits_{\Omega_r} {  C(x) P_t^2  dx} =
\tfrac{\|C(x)\|_1}{C_{ld}} a(t).
\end{equation}
Due to (\ref{as0}) and (\ref{c-000}), from (\ref{c-3-2})   we get
\begin{equation}\label{c-4}
 y'(t) \leqslant C_2 y(t)(C_1 y^{\tfrac{1}{2}}(t) - C_3)
\end{equation}
provided $y(t_1) < (\frac{C_3}{C_1})^2$, where $C_3 = 1 - \frac{\|C(x)\|_1}{C_{ld}} > 0$. So,
$$
y_2(t) \leqslant \tfrac{\gamma^2 y_2(t_1)}{(y_2^{\tfrac{1}{2}}(t_1) +(\gamma - y_2^{\tfrac{1}{2}}(t_1) ) e^{\tfrac{\alpha \beta}{2}(t-t_1)})^2}
 \ \forall\, t \in [t_1,t_2) \text{ if } y(t_1) < \gamma^2,
$$
where $\alpha =C_1 C_2$  and $\gamma = \frac{C_3}{C_1} < \beta$.
By the periodicity $y_2(t_2) = y_1(0)$, we get
\begin{multline*}
e^{\tfrac{\alpha \beta}{2}t_2} \leqslant \tfrac{y_2^{\tfrac{1}{2}}(t_1)(\gamma - y_1^{\tfrac{1}{2}}(0))}{
 y_1^{\tfrac{1}{2}}(0)(\gamma - y_2^{\tfrac{1}{2}}(t_1)) }  e^{\tfrac{\alpha \beta}{2}t_1}
 \leqslant  \tfrac{y_2^{\tfrac{1}{2}}(t_1)(\gamma - y_1^{\tfrac{1}{2}}(0))}{
 y_1^{\tfrac{1}{2}}(0)(\gamma - y_2^{\tfrac{1}{2}}(t_1)) } \, \tfrac{y_1^{\tfrac{1}{2}}(0)(\beta - y_1^{\tfrac{1}{2}}(t_1))}{
 y_1^{\tfrac{1}{2}}(t_1)(\beta - y_1^{\tfrac{1}{2}}(0)) } =    \\
\tfrac{ \gamma - y_1^{\tfrac{1}{2}}(0) }{
 \beta - y_1^{\tfrac{1}{2}}(0) } \, \tfrac{y_2^{\tfrac{1}{2}}(t_1)(\beta - y_1^{\tfrac{1}{2}}(t_1))}{
 y_1^{\tfrac{1}{2}}(t_1)(\gamma - y_2^{\tfrac{1}{2}}(t_1)) }   \Rightarrow
t_1 < t_2 \leqslant T^*:= \tfrac{2}{\alpha \beta} \ln \Bigl[ \tfrac{y_2^{\tfrac{1}{2}}(t_1)(\beta - y_1^{\tfrac{1}{2}}(t_1))}{
\kappa y_1^{\tfrac{1}{2}}(t_1)(\gamma - y_2^{\tfrac{1}{2}}(t_1)) } \Bigr],
\end{multline*}
where $\kappa : = \frac{ \beta - y_1^{\frac{1}{2}}(0) }{ \gamma - y_1^{\frac{1}{2}}(0) } > 1 $, provided
$$
y_2(t_1) > \tfrac{\gamma^2 \kappa^2 e^{ \alpha \beta t_1} y_1(t_1)}{(\beta + (\kappa e^{\tfrac{\alpha \beta}{2}t_1} -1)y_1^{\tfrac{1}{2}}(t_1))^2},
\ y_1(0) < \gamma^2
\text{ and } \tfrac{\|C(x)\|_1}{C_{ld}} < 1.
$$
As a result,  we obtain the main a priori estimate
\begin{equation}\label{apr-0}
\tfrac{1}{2}  \int \limits_{\Omega_r} {  q(x,P) \, P_x^2 dx} +
 \iint \limits_{Q_{t}} { C(x) P_t^2  dx} \leqslant C_0 < \infty
\end{equation}
for all $t \in(0,t_2) $.

Now we show that the solution constructed above is unique. We use a proof by contradiction.

\textbf{\emph{Uniqueness}:} Let  $u$ and $v$ be two
solutions to the problem (\ref{a-1})--(\ref{a-2-2}). Let us denote by
$w =u - v$ satisfying
\begin{equation} \label{w-1}
C(x) w_t =  \Bigl(q(x,u) \, w_x +
\tfrac{1}{2} ( q_u(x,u) + q_v(x,v))\,w v_x \Bigr)_x
\end{equation}
\begin{equation}\label{w-2}
\text{with initial data} \quad  w(x,0) = 0.
\end{equation}
Multiplying (\ref{w-1}) by $w(x,t)$ and integrating along $\Omega_r$, due to
(\ref{a-2-1}), we have
\begin{multline}\label{w-3}
\tfrac{1}{2}\tfrac{d}{dt} \int \limits_{\Omega_r} { C(x) w^2(x,t) dx} + \tfrac{1}{2} \tfrac{d}{dt} [C_{ls} w^2(0,t)]
+
\int \limits_{\Omega_r} {  q(x,u) \, w_x^2  dx} =   \\
 -  \int \limits_{\Omega_r} {  \tfrac{1}{2} ( q_u(x,u) + q_v(x,v) ) v_x w w_x  dx} \leqslant \\
 \tfrac{1}{2} \Bigl(\int \limits_{\Omega_r} { q(x,u) \, w_x^2  dx}  \Bigr)^{\tfrac{1}{2}}
\mathop {\sup} \limits_{\Omega_r} | q(x,v) \, v_x |
\mathop {\sup} \limits_{\Omega_r} \left[  \tfrac{ (q_u(x,u) + q_v(x,v))^2 } {  C(x) q(x,u)  q^2(x,v) }  \right]^{\tfrac{1}{2}}
\Bigl( \int \limits_{\Omega_r} { C(x) w^2(x,t) dx}  \Bigr)^{\tfrac{1}{2}}
  \end{multline}
for all $t \in (0,t_1)$. Using Cauchy inequality, boundedness of $u$ and $v$, $q(x,v) \, v_x
\in L^2(0,t_1; H^1(\Omega_r))$, due to Gr\"{o}nwall lemma, we get
\begin{equation}\label{w-4}
\int \limits_{\Omega_r} { C(x) w^2(x,t) dx} \leqslant 0 \Rightarrow w(x,t) = 0 \ \forall\,t \in (0,t_1).
\end{equation}
On the other hand, multiplying (\ref{w-1}) by $w(x,t)$ and integrating along $\Omega_r$, due to
(\ref{a-2-2}), we have
\begin{multline}\label{w-4}
\tfrac{d}{dt} \int \limits_{\Omega_r} { C(x) w^2(x,t) dx} +
\int \limits_{\Omega_r} {  q(x,u) \, w_x^2  dx} \leqslant   \\
   \tfrac{d}{dt} [C_{ld} w^2(L_r,t)] + \tilde{C}\, \mathop {\sup} \limits_{\Omega_r} | q(x,v) \, v_x |^2
\int \limits_{\Omega_r} { C(x) w^2(x,t) dx}
  \end{multline}
for all $t \in ( t_1, t_2)$, where $\tilde{C} > 0$. As
$$
\int \limits_{\Omega_r} { C(x) w(x,t) dx} = C_{ld} w(L_r, t)  \ \forall\,t \in ( t_1, t_2),
$$
then
\begin{multline} \label{w-5}
C_{ld} w^2(L_r,t) = \tfrac{1}{C_{ld}}\Bigl(\int \limits_{\Omega_r} { C(x) w  (x,t) dx}  \Bigr)^2   \leqslant
\tfrac{\|C(x)\|_1}{C_{ld}} \int \limits_{\Omega_r} { C(x) w^2(x,t) dx}  \ \forall\,t \in ( t_1, t_2).
\end{multline}
So, integrating (\ref{w-4}) in time from $t_1$ to $t$, using (\ref{w-5}), we have
\begin{multline}\label{w-6}
(1- \tfrac{\|C(x)\|_1}{C_{ld}})  \int \limits_{\Omega_r} { C(x) w^2(x,t) dx} +  C_{ld} w^2(L_r,t_1) +
 \iint \limits_{Q_{t_1}^{t}} {  q(x,P) \, P_x^2  dx dt} \leqslant \\
\tilde{C} \int \limits_{t_1}^t { \mathop {\sup} \limits_{\Omega_r} | q(x,v) \, v_x |^2
\int \limits_{\Omega_r} { C(x) w^2(x,t) dx} dt}  \ \forall\,t \in ( t_1, t_2)
\end{multline}
provided  $ \|C(x)\|_1 <  C_{ld} $. Applying Gr\"{o}nwall lemma to (\ref{w-6}), we get
$w (x,t) = 0$ for all $t \in ( t_1, t_2)$. As a result, we obtain that $w = 0 \Leftrightarrow u = v$ for all
$t \in ( 0, t_2)$.

We proceed by obtaining well-posedness conditions for the second part of the interval, namely $(L_r,L_l)$   where, to compare to the first part $(0,L_r)$,  Neumann amd dynamical flux boundary conditions are switched.
Consider the following problem
\begin{equation}\label{b-1a}
C(x) P_t =  ( q(x,P) \, P_x)_x  
\end{equation}
with initial and boundary conditions

\begin{equation}\label{b-4}
P_1(x,0) = P_0(x)  \geqslant 0,
\end{equation}
\begin{equation}\label{b-2-1}
C_{rs}  P_{1,t} =  q(x,P_1) \, P_{1,x}  \text{ at } x = L_r,\
P_{1,x} = 0 \text{ at } x = L_l
\end{equation}
for all $t \in (0,t_1)$;

\begin{equation}\label{b-4-2}
P_2(x,t_1) = \left \{
\begin{gathered}
P_1(x,t_1) \ \forall\,x \in [L_r,L_l -\epsilon],  \\
g(x,t_1) \ \forall\,x \in (L_l -\epsilon,L_l],
\end{gathered}  \right.
\end{equation}
where $g(x,t_1)$ such that
\begin{equation}\label{b-f}
\int \limits_{\Omega_l^{\epsilon}} { C(x)
g(x,t_1) dx} - C_{rd} g(L_l,t_1) =
\int \limits_{\Omega_l^{\epsilon}} { C(x) P_1 (x,t_1) dx}
+ C_{rs} P_1(0,t_1),
\end{equation}
\begin{equation}\label{b-2-2a}
P_{2,x} = 0 \,\, \text{ at } \,\,  x = L_r, \  C_{rd} P_{2,t} =  q(x,P_2) \, P_{2,x} \,\, \text{ at } \,\, x = L_l
\end{equation}
for all $t \in (t_1,t_2)$;

\begin{equation}\label{b-4-2-00}
P_1(x,t_2) = \left \{
\begin{gathered}
\tilde{g}(x,t_2) \ \forall\,x \in [L_r, L_r+\epsilon],  \\
P_2(x,t_2) \ \forall\,x \in ( L_r + \epsilon,L_l],
\end{gathered}  \right.
\end{equation}
where $\tilde{g}(x,t_2)$ such that
\begin{equation}\label{b-f-00}
\int \limits_{L_r}^{L_r+\epsilon} { C(x)
\tilde{g}(x,t_2) dx} + C_{rs} \tilde{g}(L_r,t_2) =
\int \limits_{L_r}^{L_r+\epsilon} { C(x) P_2 (x,t_2) dx}
- C_{rd} P_2(L_l,t_2).
\end{equation}

Let
\begin{equation}\label{b-sol}
P(x,t) = \left \{
\begin{gathered}
P_1(x,t) \ \forall\,t \in (0,t_1),  \\
P_2(x,t) \ \forall\,t \in [t_1,t_2),
\end{gathered}  \right.
\end{equation}
be a solution to problem (\ref{b-1a})--(\ref{b-2-2a}) for all $t \in (0,t_2)$ such that
\begin{equation}\label{b-period}
P(x,t) = P(x,t+t_2).
\end{equation}
Introduce the following notations
$$
Q_{t_1}^{t_2} = \Omega_l \times (t_1,t_2),\  Q_{t_2} = \Omega_l \times (0,t_2),
\ \Omega_l := (L_r,L_l), \  \Omega_l^{\epsilon} := ( L_l - \epsilon, L_l).
$$

\begin{definition}\label{b-def}
A non-negative function $P(x,t)$ is said to be a periodic solution of the problem (\ref{b-1a})--(\ref{b-2-2a}),
i.\,e. $P(x,t) = P(x,t+t_2)$, if
$$
P \in C_{x,t}^{\tfrac{1}{2},\tfrac{1}{4}}(\bar{Q}_{t_2}) \cap  L^{\infty}(0,t_2;H^1(\Omega_l)) ,\ C^{\tfrac{1}{2}}(x)P_t \in L^2(Q_{t_2}),
$$
$$
q(x,P) \, P_x  \in L^2 (0,t_2;H^1(\Omega_l)),
$$
and $P(x,t)$ satisfies equation (\ref{b-1a})  in
the sense that
$$
\iint \limits_{Q_{t_2}} {C(x) P_t \phi  \,dx dt} -  \iint \limits_{Q_{t_2}}
{ \bigl( q(x,P) \, P_x \bigr)_x \phi  \,dx dt} =0
$$
for any $\phi \in L^2(Q_{t_2})$ and  $\phi(x,0) = \phi(x,t_2)$.
\end{definition}

\begin{theorem} \label{b-ex}
Assume that
$$
C(x) \in C(\bar{\Omega}_l) :  \inf C(x) > 0, \ \|C(x)\|_1   < C_{rd}; \ C_{rs} \geqslant 0,  \ C_{rd} > 0;
$$
$$
 \inf N(x) > 0, \,\,  \inf A_0(x) > 0,\  g(x,t_1) \in C(\bar{\Omega}_l^\epsilon),
$$
and initial data $P_0(x) \in H^1(\Omega_l)$ is non-negative satisfying
$$
\int \limits_{\Omega_l} { q(x,P_0) \, P_{0,x}^2 dx} \leqslant
\tfrac{1}{64 \pi \mu^2 |\Omega_l|} \tfrac{ (C_{rd} - \| C(x)\|_1) }{C^2_{rd}} \mathop {\inf} \limits_{\Omega_l}\left| \tfrac{N(x) A_0^4(x)}{ C^2(x)}\right|,
$$
$$
M:= \int \limits_{\Omega_l} { C(x) P_0 (x) dx} + C_{rs}P_0 (0)   > 0 ,
$$
then the problem (\ref{b-1a})--(\ref{b-2-2a}) admits a unique non-negative solution
in the sense of Definition~\ref{b-def}.
\end{theorem}

The proof of Theorem~\ref{b-ex} is similar to the one of Theorem~\ref{ex}.
Finally, to get well-posedness for the whole interval  $(0,L_l)$, the restrictions on the parameter values obtained in  Theorem~\ref{ex} should be combined with the restrictions obtained in  Theorem~\ref{b-ex}.

\subsection{Existence and uniqueness of the nonnegative solution in Fontan case.}
Let us introduce the following notation
$$
Q_{t_1}^{t_2} = \Omega_L \times (t_1,t_2),\  Q_{t_2} = \Omega_L \times (0,t_2),
\ \Omega_L := (0,L_l), \  \Omega_L^{\epsilon} := ( L_l - \epsilon, L_l).
$$
Consider the following equation
\begin{equation}\label{cc-1}
C(x) P_t =  \bigl( q(x,P) \, P_x \bigr)_x  
\end{equation}
with initial and boundary conditions

\begin{equation}\label{cc-4}
P_1(x,0) = P_0(x)  \geqslant 0,
\end{equation}
\begin{equation}\label{cc-2-1}
C_{ls}  P_{1,t} =   q(x,P_1) \, P_{1,x}  \text{ at } x = 0,\
P_{1,x} = 0 \text{ at } x = L_l
\end{equation}
for all $t \in (0,t_1)$;

\begin{equation}\label{cc-4-2}
P_2(x,t_1) = \left \{
\begin{gathered}
P_1(x,t_1) \ \forall\,x \in [0,L_l -\epsilon],  \\
k(x,t_1) \ \forall\,x \in (L_l -\epsilon,L_l],
\end{gathered}  \right.
\end{equation}
where $k(x,t_1)$ such that
\begin{equation}\label{cc-f}
\int \limits_{\Omega_L^{\epsilon}} { C(x)
k(x,t_1) dx} - C_{ld} k(L_l,t_1) =
\int \limits_{\Omega_L^{\epsilon}} { C(x) P_1 (x,t_1) dx}
+ C_{ls} P_1(0,t_1),
\end{equation}
\begin{equation}\label{cc-2-2}
P_{2,x} = 0 \,\, \text{ at } \,\,  x = 0, \  C_{ld} P_{2,t} =   q(x,P_2) \, P_{2,x} \,\,  \text{ at } \,\, x = L_l
\end{equation}
for all $t \in (t_1,t_2)$;

\begin{equation}\label{cc-4-2-00}
P_1(x,t_2) = \left \{
\begin{gathered}
\tilde{k}(x,t_2) \ \forall\,x \in [0,  \epsilon],  \\
P_2(x,t_2) \ \forall\,x \in ( \epsilon,L_l],
\end{gathered}  \right.
\end{equation}
where $\tilde{k}(x,t_2)$ such that
\begin{equation}\label{cc-f-00}
\int \limits_{0}^{ \epsilon} { C(x)
\tilde{k}(x,t_2) dx} + C_{ls} \tilde{k}(0,t_2) =
\int \limits_{0}^{\epsilon} { C(x) P_2 (x,t_2) dx}
- C_{ld} P_2(L_l,t_2).
\end{equation}

Let
\begin{equation}\label{cc-sol}
P(x,t) = \left \{
\begin{gathered}
P_1(x,t) \ \forall\,t \in (0,t_1),  \\
P_2(x,t) \ \forall\,t \in [t_1,t_2),
\end{gathered}  \right.
\end{equation}
be a solution to problem (\ref{cc-1})--(\ref{cc-2-2}) for all $t \in (0,t_2)$ such that
\begin{equation}\label{cc-period}
P(x,t) = P(x,t+t_2).
\end{equation}

\begin{definition}\label{cc-def}
A non-negative function $P(x,t)$ is said to be a periodic solution of the problem (\ref{cc-1})--(\ref{cc-2-2}),
i.\,e. $P(x,t) = P(x,t+t_2)$, if
$$
P \in C_{x,t}^{\tfrac{1}{2},\tfrac{1}{4}}(\bar{Q}_{t_2}) \cap  L^{\infty}(0,t_2;H^1(\Omega_L)) ,\ C^{\tfrac{1}{2}}(x)P_t \in L^2(Q_{t_2}),
$$
$$
q(x,P) \, P_x  \in L^2 (0,t_2;H^1(\Omega_L)),
$$
and $P(x,t)$ satisfies equation (\ref{cc-1})  in
the sense that
$$
\iint \limits_{Q_{t_2}} {C(x) P_t \phi  \,dx dt} - \iint \limits_{Q_{t_2}}
{ ( q(x,P) \, P_x)_x \phi  \,dx dt} =0
$$
for any $\phi \in L^2(Q_{t_2})$ and  $\phi(x,0) = \phi(x,t_2)$.
\end{definition}

\begin{theorem} \label{cc-ex}
Assume that
$$
C(x) \in C(\bar{\Omega}_L) :  \inf C(x) > 0, \ \|C(x)\|_1   < C_{ld}; \ C_{ls} \geqslant 0,  \ C_{ld} > 0;
$$
$$
 \inf N(x) > 0, \, \, \inf A_0(x) > 0,\  g(x,t_1) \in C(\bar{\Omega}_L^\epsilon),
$$
and initial data $P_0(x) \in H^1(\Omega_L)$ is non-negative satisfying
$$
\int \limits_{\Omega_L} {   q(x,P_0) \, P_{0,x}^2 dx} \leqslant
\tfrac{1}{64 \pi \mu^2 |\Omega_L|} \tfrac{ (C_{ld} - \| C(x)\|_1) }{C^2_{ld}} \mathop {\inf} \limits_{\Omega_L}\left| \tfrac{N(x) A_0^4(x)}{ C^2(x)}\right|,
$$
$$
M:= \int \limits_{\Omega_L} { C(x) P_0 (x) dx} + C_{ls}P_0 (0)   > 0 ,
$$
then the problem (\ref{cc-1})--(\ref{cc-2-2}) admits a unique positive solution
in the sense of Definition~\ref{cc-def}.
\end{theorem}
The proof of Theorem~\ref{cc-ex} is similar to that of Theorem~\ref{ex}.

\vspace{1cm}
\textbf{Acknowledgements}\\
{\small The authors would like to thank the organizers and participants of the Industrial Problem Solving Workshop (2016, Fields Institute) where the problem studied in this article was originally proposed, especially Rob Andrews, Almut Burchard, Itamar Halevy, Greg Lewis, Faizan Khalid Mohsin, Pinaki Mondal, and Siv Sivaloganathan. We also acknowledge Lucy Roche and Amine Mazine for their contributions to the clinical aspects of this project.
This work was partially supported by a grant from the Simons Foundation (no. 277088 to M. Chugunova) and by a grant (NSF-DMS 1515130 to J.P. Keener).}

\appendix
\section{On super- and sub-solutions}
The second order parabolic equation admits super- and sub-solutions.
Here we show how to constract them for a  class of bounded initial values $P_0(x)$.
By the comparison principle this implies that $P_{sub} \leq P \leq P_{sup}$.
First of all, we consider the problem for the sub-solution:
\begin{equation}\label{s-1}
 P_t =  d_1^2 \, P_{xx},
\end{equation}
where
$$
d_1^2 = \tfrac{k_1}{\inf C(x)},
 \ k_1 = \inf [q(x,P)],
$$
with initial and boundary conditions
\begin{equation}\label{s-2}
P_1(x,0) = \inf P_0(x)  \geqslant 0,
\end{equation}
\begin{equation}\label{s-3}
P_{1,t} =  \tfrac{k_1}{C_{ls}}  P_{1,x}  \text{ at } x = 0,\
P_{1,x} = 0 \text{ at } x = L_r
\end{equation}
for all $t \in (0,t_1)$;
\begin{equation}\label{s-5}
P_{2,x} = 0 \text{ at } x = 0, \  P_{2,t} =  \tfrac{k_1}{C_{ld} }   P_{2,x}  \text{ at } x = L_r
\end{equation}
for all $t \in (t_1,t_2)$. The problem (\ref{s-1})--(\ref{s-3}) has a particular solution
$$
P_1(x,t) = [a_0 \cos( \tfrac{\lambda_{1}}{d_1}x) +
b_0 \sin(\tfrac{\lambda_{1}}{d_1}x)] e^{- \lambda_{1}^2t},
$$
where $\lambda_{1}$ satisfies
$$
\tan ( \tfrac{\lambda_{1}}{d_1}L_r) = - \tfrac{d_1 C_{ls}}{k_1} \lambda_{1} .
$$
The problem (\ref{s-1}), 
(\ref{s-5}) has a particular solution
$$
P_2(x,t) =  c_0 \cos( \tfrac{\lambda_{2}}{d_1}x)  e^{- \lambda_{2}^2 t},
$$
where $\lambda_{2}$ is the solution of the following equation
$$
\tan ( \tfrac{\lambda_{2}}{d_1}L_r) =   \tfrac{d_1 C_{ld}}{k_1} \lambda_{2} .
$$
If $P_2(x,t_2)=P_1(x,0)$ then
$$
P_{sub}(x,t) = \inf P_0(x)   \quad   \text{for all} \quad t \in [0,t_2]
$$
provided
$$
b_0 = 0, \ \lambda_1 = \lambda_2 = 0, \,\,  a_0 =c_0 = \inf P_0(x).
$$

Next, we consider the problem for the super-solution:
\begin{equation}\label{s-6}
 P_t =  d_2^2  P_{xx},
\end{equation}
where
$$
d_2^2 = \tfrac{k_2}{\sup C(x)},
 \ k_2 = \sup [q(x, p)],
$$
with initial and boundary conditions
\begin{equation}\label{s-7}
P_1(x,0) = \sup P_0(x)  > 0,
\end{equation}
\begin{equation}\label{s-8}
P_{1,t} =  \tfrac{k_2}{C_{ls}}  P_{1,x} \,\, \text{ at }  \,\, x = 0,\
P_{1,x} = 0 \,\, \text{ at } \,\, x = L_r
\end{equation}
for all $t \in (0,t_1)$;
\begin{equation}\label{s-10}
P_{2,x} = 0 \text{ at } x = 0, \  P_{2,t} =  \tfrac{k_2}{C_{ld} }   P_{2,x}  \text{ at } x = L_r
\end{equation}
for all $t \in (t_1,t_2)$.
The problem (\ref{s-6})--(\ref{s-8}) has a particular solution
$$
P_1(x,t) = [\tilde{a}_0 \cos( \tfrac{\tilde{\lambda}_{1}}{d_2}x) +
\tilde{b}_0 \sin(\tfrac{\tilde{\lambda}_{1}}{d_2}x)] e^{- \tilde{\lambda}_{1}^2t},
$$
where $\tilde{\lambda}_{1}$ satisfies
$$
\tan ( \tfrac{\tilde{\lambda}_{1}}{d_2}L_r) = - \tfrac{d_2 C_{ls}}{k_2} \tilde{\lambda}_{1} .
$$
The problem (\ref{s-6}), 
(\ref{s-10}) has a particular solution
$$
P_2(x,t) =  \tilde{c}_0 \cos( \tfrac{\tilde{\lambda}_{2}}{d_2}x)  e^{- \tilde{\lambda}_{2}^2 t},
$$
where $\tilde{\lambda}_{2}$ satisfies
$$
\tan ( \tfrac{\tilde{\lambda}_{2}}{d_2}L_r) =   \tfrac{d_2 C_{ld}}{k_2} \tilde{\lambda}_{2} .
$$
Now we need to impose conditions for periodicity of sub-solutions.
If $P_2(x,t_2)=P_1(x,0)$ then
$$
P_{sup}(x,t) = \sup P_0(x) \quad  \text{for all} \quad t \in [0,t_2]
$$
provided
$$
\tilde{b}_0 = 0, \,\, \tilde{\lambda}_1 = \tilde{\lambda}_2 = 0, \,\,  \tilde{a}_0 =\tilde{c}_0 = \sup P_0(x).
$$

\bibliographystyle{alpha}
\bibliography{fontan}

\end{document}